\documentclass{amsart}
\usepackage{amscd, amssymb, amsmath}
\usepackage[arrow,matrix,curve,cmtip,ps]{xy}

\newcommand{\Lmin}{\Lambda^{\min}}

\def\bp(#1;#2,#3){{\sigma^{#2}({#1})}}

\newcommand{\NN}{\mathbb N}
\newcommand{\TT}{\mathbb T}
\newcommand{\ZZ}{\mathbb Z}

\newcommand{\Ee}{{\mathcal E}}
\newcommand{\Gg}{{\mathcal G}}
\newcommand{\Tt}{{\mathcal T}}

\DeclareMathOperator{\Aut}{Aut}
\DeclareMathOperator{\lsp}{span}
\DeclareMathOperator{\Mor}{Mor}
\DeclareMathOperator{\Obj}{Obj}
\DeclareMathOperator{\Ext}{Ext}
\DeclareMathOperator{\FE}{{\mathcal{FE}}}
\DeclareMathOperator{\supp}{supp}

\theoremstyle{plain}
\newtheorem{theorem}{Theorem}[section]
\newtheorem{cor}[theorem]{Corollary}
\newtheorem{lemma}[theorem]{Lemma}
\newtheorem{prop}[theorem]{Proposition}
\theoremstyle{remark}
\newtheorem{rmk}[theorem]{Remark}
\newtheorem{rmks}[theorem]{Remarks}

\newtheorem{notation}[theorem]{Notation}
\theoremstyle{definition}
\newtheorem{dfn}[theorem]{Definition}

\newtheorem{example}[theorem]{Example}
\newtheorem{examples}[theorem]{Examples}

\numberwithin{equation}{section}

\begin{document}

\title[{\boldmath Higher-Rank Graph $C^*$-algebras}]{{\boldmath Higher-Rank Graph $C^*$-algebras: An Inverse Semigroup
and Groupoid Approach}}
\author{Cynthia Farthing}
\author{Paul S. Muhly}
\author{Trent Yeend}
\address{Department of Mathematics \\ University of Iowa \\ IA 52242 \\ U.S.A.}
\address{Department of Mathematics \\ University of Newcastle \\  NSW  2308 \\ AUSTRALIA}

\email{cfarthin@math.uiowa.edu} \email{pmuhly@math.uiowa.edu}
\email{trent.yeend@newcastle.edu.au}

\thanks{The research of the first two authors was supported in part by a grant from the National Science Foundation,
DMS-0070405.}

\date{\today}


\begin{abstract}
We provide inverse semigroup and groupoid models for the Toeplitz and Cuntz-Krieger
algebras of \emph{finitely aligned} higher-rank graphs. Using these models, we prove a
uniqueness theorem for the Cuntz-Krieger algebra.
\end{abstract}

\maketitle

\section{Introduction}

A higher-rank graph is a countable category $\Lambda$ endowed with a degree functor $d:\Lambda\rightarrow\mathbb{N}^{k}$
satisfying the \emph{unique factorization property}: For all $\lambda\in\Lambda$ and $m,n\in\NN^k$ with
$d(\lambda)=m+n$, there are unique elements $\mu,\nu\in \Lambda$ such that $d(\mu)=m$, $d(\nu)=n$
and $\lambda=\mu\nu$. The rank of $\Lambda$ is $k$ and for this
reason, $(\Lambda,d)$ is also called a $k$\emph{-graph}. A $1$-graph
is simply the finite-path category generated freely by an ordinary directed
graph.\footnote{For the purpose of motivation, we discuss the theory of $1$-graphs at some length in the next section.
For a very readable account of graph $C^*$-algebras, including the rudiments of the $C^*$-algebras of $k$-graphs, we
recommend the CBMS lectures by Iain Raeburn \cite{iR05}.} In \cite{KP}, Kumjian and Pask introduced the notion of a
higher-rank graph in order to capture the essential features of the
$C^{\ast}$-algebras that Robertson and Steger
associated to buildings \cite{RS96,RS99,RS01,RS03} and to provide links between these and higher
order shift dynamical systems. (See \cite{KP03} also.)

The $C^{\ast}$-algebras associated to higher-rank graphs are generalizations
of ordinary graph $C^{\ast}$-algebras in that they are generated
by families of partial isometries $\{ s_{\lambda}\}_{\lambda\in\Lambda}$
that satisfy certain relations that have received a lot of attention
in recent years. Kumjian and Pask defined and studied the $C^{\ast}$-algebra
of $(\Lambda,d)$, $C^{\ast}(\Lambda)$, in terms of a certain type
of groupoid that encodes the graph. They were motivated by, and generalized,
the theory in \cite{KPRR} which has been the source of considerable
inspiration in our subject. However, just as in the setting of
ordinary graphs, where the groupoid techniques of \cite{KPRR} require
hypotheses that rule out many interesting examples, the work of Kumjian
and Pask requires hypotheses that place significant limitations on
the nature of the $k$-graphs that may be analyzed. Our first objective
in this paper, then, is to overcome the limitations that Kumjian and
Pask place on their $k$-graphs and to show how to build a groupoid
that gives the $C^{\ast}$-algebra of an arbitrary $k$-graph subject
only to the condition that it is ``finitely aligned" (see
Definition~\ref{dfn:finitely aligned}). This condition seems
to lie at the natural ``boundary" of the subject. That
is, with or without the use of groupoids, little can be said about
$k$-graphs that are not finitely aligned.

We were motivated in part by the important contribution of Paterson
\cite{P2} in which he successfully circumvented the limitations of
\cite{KPRR} that involve finiteness hypotheses on the graphs under
consideration by first introducing an inverse semigroup that is naturally
attached to the graph. Once this inverse semigroup is identified,
he constructed a groupoid from it using technology that he and others
have developed and which he exposed thoroughly in \cite{P1}. Paterson's
success inspired our approach here. Given a $k$-graph $(\Lambda,d)$,
we first build a natural inverse semigroup from $\Lambda$, $S_{\Lambda}$.
However, in contrast to the rank-$1$ setting of \cite{P2}, the
groupoid we want is fairly far removed from the universal groupoid of $S_{\Lambda}$.
Rather, it is obtained directly from a natural action of $S_{\Lambda}$ on
a certain ``infinite-path space" and is realized
in terms of the \emph{sheaf of germs} of the action (cf. \cite{AR}).

One extra benefit of our analysis is that we obtain a groupoid
presentation of the \emph{Toeplitz algebra} of $(\Lambda,d)$, $\Tt C^{*}(\Lambda)$.  Another is that we overcome the
limitation of \cite{P2} whereby the
graphs considered are free of sources. That is, our theory gives
an extension of Paterson's analysis even when restricted to the context of
ordinary graphs.

We note, too, that under the hypotheses invoked by Kumjian and Pask \cite{KP}, our analysis is different from theirs.
They build a groupoid directly from the $k$-graph; we obtain a groupoid by first considering an inverse semigroup.
Their groupoid and ours are the same, however, under their hypotheses. (See Remarks \ref{rmk:restricted groupoid
notes}.)

The paper is arranged as follows.  In the next section, we discuss some of the features of $1$-graphs that inspired
our analysis. The discussion here is informal and incomplete. Detailed work begins in Section~\ref{sec:preliminaries},
where some of the basic facts about $k$-graphs are exposed, notation is
set up, and basic facts about the $C^{*}$-algebras we will study
are presented. Also, two propositions regarding the structure of $k$-graphs
are presented to be used in Section~\ref{sec:the action}. As will
be seen in due course, they lie at the heart of our analysis.

In Section~\ref{sec:inverse semigroup}, we restrict our attention
to finitely aligned $k$-graphs $(\Lambda,d)$ (Definition~\ref{dfn:finitely aligned}),
and define our inverse semigroup $S_{\Lambda}$ from the structure
of paths in $\Lambda$. We then define a second countable,
locally compact, Hausdorff space $X_{\Lambda}$ and an inverse semigroup
action $\theta$ of $S_{\Lambda}$ on $X_{\Lambda}$. The space $X_{\Lambda}$
comprises all paths constructed on $\Lambda$ --- finite, infinite and partially
infinite --- and $S_{\Lambda}$ acts naturally by ``removing and adding
initial segments''.

Section~\ref{sec:the groupoid} defines the groupoid of germs (cf.
\cite{AR}) $\Gg_{\Lambda}$ of the system $(X_{\Lambda},S_{\Lambda},\theta)$.
With a topology naturally arising from the action and topology on
$X_{\Lambda}$, $\Gg_{\Lambda}$ becomes an ample groupoid with unit
space $X_{\Lambda}$, and $S_{\Lambda}$ is identified with an inverse
semigroup of ample subsets of $\Gg_{\Lambda}$. It is then shown that
$C^{*}(\Gg_{\Lambda})$ is isomorphic to the Toeplitz algebra $\Tt C^{*}(\Lambda)$.
A closed invariant subset $\partial\Lambda$ of the unit space $\Gg_{\Lambda}^{(0)}$
is then identified and it is shown that $C^{*}(\Gg_{\Lambda}|_{\partial\Lambda})$
is isomorphic to the $C^{*}$-algebra $C^{*}(\Lambda)$.

Finally, in section~\ref{sec:CK uniqueness} we present a Cuntz-Krieger
uniqueness theorem for finitely aligned $k$-graphs. A similar theorem
was given in \cite{RSY2}, and we compare the two results here.

\section{Motivation}
\label{sec:Motivation}

In this section we call attention to the salient features of graph $C^*$-algebras that are the inspiration for the
present work.  Our intention is to provide an outline of the crucial points of our analysis, especially to help those
unfamiliar with groupoid and inverse semigroup methods in operator algebra.  Those familiar with the theory of graph
$C^*$-algebras, and with the associated theories of inverse semigroups and groupoids, may skip directly to the next
section.

Let $E = (E^0,E^1,r,s)$ be an ordinary directed graph.  This means that $E^0$ and $E^1$ are sets, called respectively
the set of
\emph{vertices} and the set of \emph{edges}, and that $r$ and $s$ are functions from $E^1$ to $E^0$, called
respectively the \emph{range} and \emph{source} maps.  We will assume here that our graphs have countable vertex and
edge sets. From $E$ we build the associated \emph{finite-path category} $E^*$.  This is just the free category
generated by $E$. It may and will be viewed as the collection of finite words $(e_1,e_2, \dots ,e_n)$ over $E^1$,
where $s(e_i)=r(e_{i+1})$ for all $i \le n-1$, together with the vertices $v \in E^0$.  The range and source maps
extend in the obvious way to $E^*$ so that $E^*$ is also a graph.  Composition in $E^*$ is defined through
concatenation. Thus, paths $(e_1,e_2, \ldots ,e_n)$ and $(f_1,f_2, \ldots ,f_m)$ are composable if and only if
$s(e_n)=r(f_1)$, and in this event, their product is $(e_1,e_2, \ldots ,e_n, f_1,f_2, \ldots ,f_m)$; if $v$ is a
vertex and $e$ is of the form $(e_1,e_2, \ldots ,e_n)$, then $v$ and $e$ are composable if and only if $v = r(e)$, and
in this event, $ve$ is just $e$.

The category $E^*$ is countable and is endowed with a \emph{degree functor} $d$ from $E^*$ to the semigroup, or small
category, of non-negative integers, $\NN$. Namely, $d$ describes the length of paths in $E^*$, so $d(v) = 0$ for all
vertices $v$, and $d(e)=n$ for $e$ of the form $(e_1,\dots,e_n)$. Further, the ``freeness'' of the construction of
$E^*$ is expressed by the fact that if $d(e) = n_1 + n_2$, $n_1,n_2 \in \NN$, then there are \emph{unique} paths $e_1$
and $e_2$ in $E^*$ such that $d(e_i) = n_i$, $i = 1, 2$, and $e = e_1e_2$.  Thus $E^*$ and $d$ satisfy the
\emph{unique factorization property} that we mentioned at the outset, and $(E^*,d)$ is a $1$-graph.  The
innocuous-seeming factorization property is essential for the role graphs play in operator algebra.

Given a graph $E$, one would like to build a $C^*$-algebra $C^*(E)$ that is generated by a family $\{s_e\}_{e\in E^*}$
which, at the very least, consists of partial isometries (i.e., $s_e^*s_e$ is a projection for each $e$) satisfying
$s_e s_f = s_{ef}$ for all composable $e,f \in E^*$.  That is, one would like to build a $C^*$-algebra that codifies
the representation theory of $E^*$ by partial isometries.  Experimental investigation reveals that the crucial
elements in $C^*(E)$ are or ought to be $s_es_f^*$, $e,f \in E^*$, and that one would like the products $s_es_f^*$ to
behave like matrix units when $e$ and $f$ have the same degree.  Further, in order for the matrix units  $s_es_f^*$ of
one degree to be linked nicely with the matrix units of a different degree, one is led naturally to require the
partial isometries $s_e$ satisfy the so-called Cuntz-Krieger condition:
\begin{equation}
s_v = \sum_{\{e \in E^1: r(e) = v\}} s_e s_e^* \label{eqn: Cuntz-Krieger}
\end{equation}
Note that since $s_v$ must be an idempotent partial isometry, it is a projection, and customarily, one
writes $p_v$ instead of $s_v$ to highlight this.  Note, too, that for equation (\ref{eqn: Cuntz-Krieger}) to make
sense, one requires that $v$ must be the range of at least one edge, so, as one says, $v$ is not a \emph{source}.
Also, since an infinite sum of projections cannot converge in a $C^*$-algebra, one requires also that the sum is
finite; i.e., one requires that $v$ is not an \emph{infinite receiver}.  Problems with sources and infinite receivers
play an important role in the subject, and this paper contributes to their solution, but for this discussion, we will
assume that our graph has no sources and that there are no infinite receivers.

The history and theory of graph $C^*$-algebras is fairly complex and involved, but to keep matters short, we jump to
Paterson's wonderful insight \cite{P2} that it is very helpful to embed $E^*$ into a certain inverse semigroup
$S_{E}$, and then to use the theory he and others have been developing to realize the $C^*$-algebra $C^*(E)$ as the
$C^*$-algebra universal for particular representations of $S_{E}$ (see \cite[Theorem~2]{P2}). In fact, in \cite{P1}
Paterson advocates that the $C^*$-algebras universal for particular representations of an inverse semigroup may be
effectively studied as the $C^*$-algebras of groupoids naturally associated to the semigroup. This is the tack we take
here. Any inverse semigroup acts by partially defined homeomorphisms on the semicharacter space of its idempotent
subsemigroup, and the collection of germs of these maps forms a groupoid. In our setting, there is a closed subset of
the semicharacter space which is invariant under the action, and so we consider the groupoid of germs associated to
the restricted action.

In the context of our graph $E$, the inverse semigroup $S_{E}$ consists of a zero element $z$ together with all pairs
$(\alpha ,\beta) \in E^* \times E^*$ such that $s(\alpha) = s(\beta)$. Multiplication
in $S_E$ is given by the formula
\begin{equation}
(\mu,\nu)(\alpha,\beta) =
\begin{cases}
(\mu \alpha',\beta) &\text{if $\alpha = \nu \alpha'$} \\
(\mu, \beta \nu') &\text{if $\nu = \alpha \nu'$} \\
z &\text{otherwise,}
\end{cases}\label{eqn:Inverse Semigroup}
\end{equation}
and involution is given by $z^*=z$ and $(e,f)^*=(f,e)$. The path category $E^*$ is embedded in $S_{E}$ via the formula $e \mapsto (e, s(e))$.

As we just mentioned, we want to think of $S_{E}$ acting by partially defined homeomorphisms on
the semicharacter space $\hat{\Ee}$ of the idempotent semigroup $\Ee$ of $S_{E}$. Clearly, $\Ee$ consists of $z$
together with all the pairs $(e,e)$, $e \in E^*$.  On the other hand, by definition, $\hat{\Ee}$ consists of all
(nonzero) semigroup homomorphisms from $\Ee$ to the multiplicative semigroup $\{0,1\}$.  With respect to the topology
of pointwise convergence, $\hat{\Ee}$ is a locally compact Hausdorff space. The semicharacter space $\hat{\Ee}$ may be
identified as the disjoint union of $\{z\}$, $E^*$ and the \emph{infinite-path space} $E^\infty$; that is,
$\hat{\Ee}$ comprises $z$, the finite paths $(e_1,e_2,\dots,e_n)$ of $E^*$, and the infinite sequences
$(e_1,e_2,\dots)$ satisfying $s(e_i)=r(e_{i+1})$ for $i\ge 1$. The infinite-path space is a closed subset of
$\hat{\Ee}$ which is invariant under the action of $S_E$. Under the identification of $\hat{\Ee}$ with $\{z\}\cup
E^*\cup E^\infty$, the restricted action of $S_E$ on $E^\infty$ is given by removing and affixing initial segments of
infinite paths according to the formula
\[
(\alpha,\beta)\cdot x = \alpha y, \quad\text{where } x=\beta y.
\]
The groupoid $\Gg_E$ formed from the germs of the action may be realized as all triples from $E^\infty\times\ZZ\times
E^\infty$ of the form $(\alpha x, d( \alpha) - d(\beta), \beta x)$, where $\alpha$ and $\beta$ are finite paths and
$x$ is an infinite path such that $s(\alpha) = s(\beta) = r(x)$. The $C^*$-algebra of $E$ may then be realized as the
$C^*$-algebra of $\Gg_{E}$.

As we indicated at the outset, a $k$-graph is a countable category $\Lambda$, endowed with a degree functor $d$ taking
values in the semigroup $\NN^k$, that satisfies the unique factorization property: for every path $\lambda\in \Lambda$
and $m,n\in\NN^k$ with $d(\lambda) = m + n$, there exist unique paths $\mu,\nu\in\Lambda$ such that $\lambda=\mu\nu$,
$d(\mu) = m$ and $d(\nu) = n$. We want to build an inverse semigroup and groupoid for each $k$-graph. If one studies
the outline of how Paterson does this for $1$-graphs, one realizes that there are two significant hurdles that are not
anticipated by his work which must be overcome.

The first hurdle arises when one tries to put an inverse semigroup structure on the collection of pairs
$(\lambda, \mu) \in \Lambda \times \Lambda$ such that $s(\lambda) = s(\mu)$. More precisely, the multiplication
structure described by equation \eqref{eqn:Inverse Semigroup} is not appropriate for general $k$-graphs with
$k \ge 2$. Our way around this problem is to look at \emph{finite sets} of such pairs satisfying certain conditions to
be described in Section~\ref{sec:inverse semigroup}. In a natural way, this collection of finite sets turns out to be
an inverse semigroup $S_{\Lambda}$ containing a natural embedding of $\Lambda$. Indeed, the image of $\Lambda$
generates $S_\Lambda$ as an inverse semigroup with a partially defined addition structure. Further, when $\Lambda=E^*$
is a $1$-graph, $S_{\Lambda}$ contains $S_{E}$ as an inverse subsemigroup, and there is a one-to-one correspondence
between representations of $S_E$ and additive representations\footnote{These are the representations of $S_\Lambda$
which preserve its partially defined addition structure; cf. \cite[page~193]{P1}.} of $S_\Lambda$.

The second hurdle concerns the semicharacter space of the subsemigroup $E(S_{\Lambda})$ of idempotents in
$S_{\Lambda}$: it appears to be much larger than the space on which we want $S_\Lambda$ to act. We pass instead
directly to an analogue of the infinite-path space. More accurately, we consider two analogues.  The first, denoted
$X_{\Lambda}$, is, in the directed graph setting, an analogue of $E^* \cup E^{\infty}$. We show that the inverse
semigroup $S_{\Lambda}$ acts by partially defined homeomorphisms on $X_{\Lambda}$.  Further, we show in
Section~\ref{sec:the groupoid} that the groupoid of germs of this action, $\Gg_{\Lambda}$, parameterizes the
\emph{Toeplitz} algebra of $\Lambda$, $\Tt C^*(\Lambda)$; this turns out to be new even in the $1$-graph setting.  The
second analogue of $E^{\infty}$ is the subset of $X_{\Lambda}$ called the space of \emph{boundary} paths, denoted
$\partial \Lambda$.  In the setting of a graph $E$ with no sources and no infinite receivers, $\partial \Lambda$ is
$E^{\infty}$.  In general, $\partial \Lambda$ is a closed subset of $X_{\Lambda}$ that is invariant under the action
of $\Gg_{\Lambda}$.  The reduction of $\Gg_{\Lambda}$ to $\partial \Lambda$, $\Gg_{\Lambda}\vert_{\partial \Lambda}$,
is our choice for the groupoid that parameterizes the \emph{Cuntz-Krieger} algebra of $\Lambda$, $C^*(\Lambda)$.

\section{Preliminaries}

\label{sec:preliminaries}

\subsection{Higher-Rank Graphs and their $C^{*}$-algebras}

Throughout the remainder of this paper, $\Lambda$ will denote a fixed $k$-graph and $d:\Lambda \rightarrow \NN^k$ will
be the associated degree functor.  As we have spelled out above, the crucial property of $d$ is the unique
factorization property: If $d(\lambda) = m + n$ in $\NN^k$, then there are unique $\mu$ and $\nu$ in $\Lambda$ such
that $d(\mu) = m$, $d(\nu)=n$, and $\lambda = \mu \nu$. Ordinarily a category $\Lambda$ is viewed as a system
$(\Obj(\Lambda),\Mor(\Lambda),r,s)$
where $\Obj(\Lambda)$ and $\Mor(\Lambda)$ are separate sets and
$r$ and $s$ are maps from the second set to the first. However,
it is convenient here to take the ``arrows only" approach to categories \cite{sM71}. So all elements of $\Lambda$
are morphisms and $\Obj(\Lambda)$ is distinguished by virtue of being
the idempotent morphisms. This perspective is especially appropriate
because of our use of the degree functor. The unique factorization property allows us to identify $\Obj(\Lambda)$ with
the
elements of $\Lambda$ that have degree zero. Because of our desire
to generalize $1$-graphs, we will also call the elements of $\Lambda$
(finite) \emph{paths}. For $m\in\NN^{k}$, we define $\Lambda^{m}:=d^{-1}(\{ m\})$,
so as we just mentioned, $\Obj(\Lambda)=\Lambda^{0}$.

We follow \cite[Section~2]{RSY1} for the basic facts about $k$-graphs
that we shall use.

\begin{notation}
\label{N:basic}
\begin{enumerate}
\item For $v\in\Lambda^{0}$ define $v\Lambda:=r^{-1}(v)$ and $\Lambda v:=s^{-1}(v)$, and for $n\in\NN^k$ define
$v\Lambda^n := \Lambda^n\cap v\Lambda$. For $\lambda,\mu\in\Lambda$ define
\[
\Lmin(\lambda,\mu):=\{(\alpha,\beta)\in\Lambda\times\Lambda:\lambda\alpha=
\mu\beta~\text{and}~d(\lambda\alpha)=d(\lambda)\vee d(\mu)\}.
\]
Here, and throughout, given $m,n\in\NN^{k}$ we write $m\vee n$
for the coordinate-wise maximum of $m$ and $n$. That is, the $i^{\text{th}}$
coordinate of $m\vee n$ is the maximum of the $i^{\text{th}}$ coordinates
of $m$ and $n$.
\item $\Lambda*_{s}\Lambda=\{(\lambda,\mu)\in\Lambda\times\Lambda:s(\lambda)=s(\mu)\}$.
\item Let $\lambda\in\Lambda$ and let $m$ and $n$ satisfy the
inequality $0\leq m\leq n\leq d(\lambda)$. Then the unique factorization
property guarantees that there are unique paths ${\lambda}_{i}$,
$i=1,2,3$, such that $d({\lambda}_{1})=m$, $d({\lambda}_{2})=n-m$, $d({\lambda}_{3})=d(\lambda)-n$ and $\lambda =
\lambda_1\lambda_2\lambda_3$. We shall write ${\lambda}(0,m)$
for ${\lambda}_{1}$, ${\lambda}(m,n-m)$ for ${\lambda}_{2}$ and
${\lambda}(n,d(\lambda))$ for ${\lambda}_{3}$.
\end{enumerate}
\end{notation}

\begin{example}\label{k-graph examples}
For $m\in(\NN\cup\{\infty\})^{k}$, define $\Omega_{k,m}$ to be the
$k$-graph with \[
\Obj(\Omega_{k,m})=\{ p\in\NN^{k}:p\leq m\},\]
 \[
\Mor(\Omega_{k,m})=\{(p,q)\in\Obj(\Omega_{k,m})\times\Obj(\Omega_{k,m}):p\leq q\},\]
 \[
r(p,q)=p,\quad s(p,q)=q,\quad d(p,q)=q-p.\]  Drawn below are
$\Omega_{2,(\infty,\infty)}$ and $\Omega_{2,(1,2)}$. In the diagrams,
edges of degree $(1,0)$ are solid; edges of degree $(0,1)$ are dashed.  In
each diagram $\lambda=((0,2),(1,2))$ and $\mu=((0,0),(0,1))$.

$$
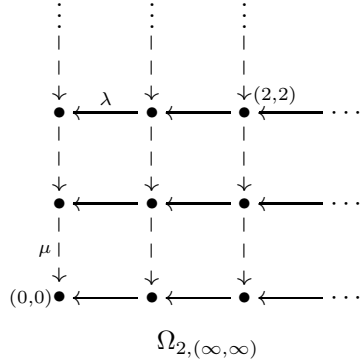
\xymatrix{\vdots\ar@{-->}[d]&\vdots\ar@{-->}[d]&\vdots\ar@{-->}[d]\\
\bullet\ar@{-->}[d]&\bullet\ar@{-->}[d]\ar[l]_{\lambda}&\bullet\ar@{-->}[d]\ar[l]&\cdots\ar[l]_>>>{(2,2)}\\
\bullet\ar@{-->}[d]_{\mu}&\bullet\ar@{-->}[d]\ar[l]&\bullet\ar@{-->}[d]\ar[l]&\cdots\ar[l]\\
\llap{$\scriptstyle (0,0)$}\bullet&\bullet\ar[l]&\bullet\ar[l]&\cdots\ar[l]}
$$
\hspace{5.75cm}$\Omega_{2,(\infty,\infty)}$

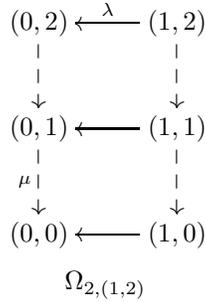
$$
\xymatrix{(0,2)\ar@{-->}[d]&(1,2)\ar[l]_{\lambda}\ar@{-->}[d]\\
(0,1)\ar@{-->}[d]_{\mu}&(1,1)\ar[l]\ar@{-->}[d]\\
(0,0)&(1,0)\ar[l]}$$
\hspace{5.75cm}$\Omega_{2,(1,2)}$

\end{example}

\begin{dfn}\label{dfn:morphism}A morphism between two \(k\)-graphs \((\Lambda_1,d_1)\) and \((\Lambda_2,d_2)\) is a
functor \(f:\Lambda_1\to\Lambda_2\) satisfying \(d_2(f(\lambda))=d_1(\lambda)\) for all \(\lambda\in\Lambda_1\).
\end{dfn}

\begin{dfn}\label{dfn:finitely aligned}
A \(k\)-graph \((\Lambda,d)\) is \emph{finitely aligned} if \(\Lmin(\lambda,\mu)\) is at most finite for all
\(\lambda,\mu\in\Lambda\).
\end{dfn}

\begin{rmk}\label{rmk:finitely aligned}
If $\Lambda$ is a 1-graph, then $\Lambda$ is automatically finitely
aligned: for $\lambda,\mu\in\Lambda$, $\Lmin(\lambda,\mu)$ is either empty
or a singleton.  Also every row-finite $k$-graph is
finitely aligned, where a $k$-graph is called row-finite if $v\Lambda^{n}$
is finite for each $v\in\Lambda^0$ and $n\in\NN^k$. \end{rmk}

\begin{dfn}\label{dfn:finite exhaustive}
Let $\Lambda$ be a $k$-graph and let $v\in\Lambda^0$.  We say a subset $E\subseteq v\Lambda$ is \emph{exhaustive} if
for every $\mu\in v\Lambda$ there exists a $\lambda\in E$ such that $\Lmin(\lambda,\mu)\not =\emptyset$.  We denote the
set of all \emph{finite exhaustive subsets} of $\Lambda$ by $\FE(\Lambda)$, and for $v\in \Lambda^0$, we define
$v\FE(\Lambda):=\{E\in \FE(\Lambda):E\subseteq v\Lambda\}$.
\end{dfn}

\begin{examples}\label{ex:finitely aligned}
\noindent
\begin{enumerate}
\item For all $m\in(\NN\cup\{\infty\})^k$ and $v\in\Omega_{k,m}^0$, any nonempty finite subset of $v\Omega_{k,m}$ is
finite exhaustive.
\item Consider the $k$-graph $\Lambda$ below:

$$
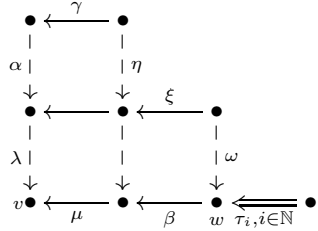
\xymatrix{\bullet\ar@{-->}[d]_{\alpha}  & \bullet\ar[l]_{\gamma}\ar@{-->}[d]^{\eta}\\
\bullet \ar@{-->}[d]_{\lambda}&\bullet \ar[l]\ar@{-->}[d]&\bullet \ar[l]_{\xi}\ar@{-->}[d]^{\omega}\\
\llap{$\scriptstyle v$}\bullet &\bullet \ar[l]^{\mu}&\smash{\underset{w}{\bullet}}\vphantom{\bullet}
\ar[l]^{\beta}&\bullet\ar@{=>}[l]^{\tau_i, i\in\NN}}$$

\noindent Dashed edges represent edges of degree $(0,1)$ and solid edges represent edges of degree $(1,0)$.  The edges
$\tau_i$ where $i\in\NN$ each have degree $(1,0)$. Any finite exhaustive subset of $w\Lambda$ must contain $w$.  The
set $\{\mu\}$ is a finite exhaustive subset of $v\Lambda$, whereas $\{\lambda\}$ is not because
$\Lmin(\lambda,\mu\beta\tau_i)=\emptyset$ for any $i\in\NN$.
\end{enumerate}
\end{examples}

\begin{dfn}\label{C-K relations}
Let \(\Lambda\) be a finitely aligned \(k\)-graph. A \emph{Toeplitz-Cuntz-Krieger \(\Lambda\)-family} in a
\(C^*\)-algebra \(B\) is a collection \(\{ t_\lambda : \lambda\in\Lambda\}\) of partial isometries in \(B\) satisfying
\begin{itemize}
\item[$(1)$]
\(\{ t_v  :  v\in\Lambda^0\}\) consists of mutually orthogonal projections;
\item[$(2)$]
\(t_\lambda t_\mu = t_{\lambda\mu}\) whenever \(s(\lambda)=r(\mu)\); and
\item[$(3)$]
\(t^*_\lambda t_\mu = \sum_{(\alpha,\beta)\in\Lmin(\lambda,\mu)}t_\alpha t^*_\beta\) for all
\(\lambda,\mu\in\Lambda\).
\end{itemize}
A \emph{Cuntz-Krieger \(\Lambda\)-family} is a Toeplitz-Cuntz-Krieger \(\Lambda\)-family \(\{t_\lambda  :
\lambda\in\Lambda\}\) which satisfies
\begin{itemize}
\item[(CK)]
\(\prod_{\lambda\in E} (t_v - t_\lambda t^*_\lambda) = 0\) for every \(v\in\Lambda^0\) and \(E\in v\FE(\Lambda)\).
\end{itemize}
\end{dfn}

Of course, the hypothesis that $({\Lambda},d)$ is finitely aligned
guarantees that the sums in Definition~\ref{C-K relations} are finite
sums, and hence make sense in any $C^{*}$-algebra. The following
remark summarizes the ontological properties of these relations and
the $C^{*}$-algebras that codify them.

\begin{rmk}\label{rmk:universal properties}
Let \((\Lambda,d)\) be a finitely aligned \(k\)-graph. Then there is a \(C^*\)-algebra \(\Tt C^*(\Lambda)\), called the
\emph{Toeplitz algebra} of \(\Lambda\), generated by a Toeplitz-Cuntz-Krieger \(\Lambda\)-family \(\{s^\Tt_\lambda :
\lambda\in\Lambda\}\), which is universal in the sense that if \(\{t_\lambda : \lambda\in\Lambda\}\) is a
Toeplitz-Cuntz-Krieger \(\Lambda\)-family in a \(C^*\)-algebra $B$, then there exists a \(C^*\)-homomorphism \(\pi :\Tt
C^*(\Lambda)\to B\) such that \(\pi(s^\Tt_\lambda)=t_\lambda\) for all \(\lambda\in\Lambda\). The Cuntz-Krieger algebra
of \(\Lambda\) is the \(C^*\)-algebra \(C^*(\Lambda)\), generated by a Cuntz-Krieger \(\Lambda\)-family \(\{s_\lambda :
\lambda\in\Lambda\}\), which is universal for Cuntz-Krieger \(\Lambda\)-families. If we let \(I\) be the ideal in \(\Tt
C^*(\Lambda)\) generated by products \(\prod_{\lambda\in E} (s^\Tt_v - s^\Tt_\lambda (s^\Tt_\lambda)^* )\) where
\(v\in\Lambda^0\) and \(E\in v\FE(\Lambda)\), then \(C^*(\Lambda)\) can be identified with \(\Tt C^*(\Lambda)/I\).
Furthermore, \cite[Proposition~2.12]{RSY2} says that the \(s^\Tt_\lambda\) and \(s_\lambda\) are nonzero for all
\(\lambda\in\Lambda\), implying that \(\Tt C^*(\Lambda)\) and \(C^*(\Lambda)\) are nontrivial.
\end{rmk}

We want to call special attention to \cite[Appendix~A]{RSY2} for
a thorough explanation of the Cuntz-Krieger relations of finitely
aligned $k$-graphs, and to \cite[Appendix~B]{RSY2} for an account
of how the theory encompasses the Cuntz-Krieger relations described in Equation~\ref{eqn: Cuntz-Krieger}.

\subsection{Extending paths}

\label{subsec:Ext}

In the study of graph $C^{*}$-algebras, i.e., in the study of $C^{*}$-algebras
defined by $1$-graphs, problems arise if the graph has sources or
sinks. Recall that a source in a graph is a vertex $v$ that does not receive
any edges, i.e., $r^{-1}(v) = \emptyset$, while $v$ is a sink if $v$ does not emit any edges, i.e., if $s^{-1}(v) =
\emptyset$. There
is now a substantial literature on how to tackle these. In a
$k$-graph $\Lambda$, a vertex $v$ is a source if
$v\Lambda^{e_i}=\emptyset$ for some $i\in\{1,\dots,k\}$ where $e_i$ is the element in $\NN^k$ with 1 in the
$i^{\text{th}}$ coordinate and 0 in every other coordinate.  Therefore
in the setting of higher-rank graphs, the situations can be
considerably more
complicated. A vertex may
receive edges from some directions
and not from others, while emitting edges in still other directions,
but not in all. The highly ramified collection of possibilities creates
numerous difficulties. One of the achievements of our approach to
the analysis of $C^{*}$-algebras associated to higher-rank graphs
is to circumvent difficulties that sources and sinks can cause. We
do not eliminate all problems, of course. Rather, we show that sources
and sinks do not prevent one from defining and analyzing a groupoid associated to such $k$-graphs.

A key tool in our analysis are the two propositions of this subsection,
Propositions~\ref{prop:C.4 RSY2} and \ref{prop:Ext transitive}.

\begin{dfn}
Given \(\lambda \in \Lambda\) and \(E \subseteq r(\lambda)\Lambda\), write \(\Ext(\lambda;E)\) for the set
\[
\bigcup_{\mu \in E} \{\alpha\in\Lambda : (\alpha,\beta) \in \Lmin(\lambda,\mu) \text{ for some } \beta\in\Lambda\}.
\]
\end{dfn}

\begin{prop}
\cite[Lemma~C.4]{RSY2}\label{prop:C.4 RSY2} Let $(\Lambda,d)$
be a finitely aligned $k$-graph, let $v\in\Lambda^{0}$, let $\lambda\in v\Lambda$,
and suppose $E\in v\FE(\Lambda)$. Then $\Ext(\lambda;E)\in s(\lambda)\FE(\Lambda)$.
\end{prop}

\begin{prop}\label{prop:Ext transitive} If $(\Lambda,d)$ is a $k$-graph, then
for $v\in\Lambda^{0}$, $E\subseteq v\Lambda$, $\lambda_{1}\in v\Lambda$
and $\lambda_{2}\in s(\lambda_{1})\Lambda$, the following equation
holds: \[
\Ext(\lambda_{2};\Ext(\lambda_{1};E))=\Ext(\lambda_{1}\lambda_{2};E).\]

\end{prop}
\begin{proof}
Let $\alpha\in\Ext(\lambda_{2};\Ext(\lambda_{1};E))$, so $(\alpha,\beta)\in\Lmin(\lambda_{2},\xi)$
for some $\xi\in\Ext(\lambda_{1};E)$ and $\beta\in\Lambda$. Then
$\lambda_{2}\alpha=\xi\beta$ and $d(\lambda_{2}\alpha)=d(\lambda_{2})\vee d(\xi)$.
Since $\xi\in\Ext(\lambda_{1};E)$, we have $(\xi,\eta)\in\Lmin(\lambda_{1},\mu)$
for some $\mu\in E$ and $\eta\in\Lambda$. So $\lambda_{1}\xi=\mu\eta$
and $d(\lambda_{1}\xi)=d(\lambda_{1})\vee d(\mu)$. Thus \[
\lambda_{1}\lambda_{2}\alpha=\lambda_{1}\xi\beta=\mu\eta\beta\]
 and \begin{align*}
d(\lambda_{1}\lambda_{2}\alpha) & =d(\lambda_{1})+d(\lambda_{2}\alpha)\\
 & =d(\lambda_{1})+(d(\lambda_{2})\vee d(\xi))\\
 & =d(\lambda_{1}\lambda_{2})\vee d(\lambda_{1}\xi)\\
 & =d(\lambda_{1}\lambda_{2})\vee(d(\lambda_{1})\vee d(\mu))\\
 & =d(\lambda_{1}\lambda_{2})\vee d(\mu),\end{align*}
 so $(\alpha,\eta\beta)\in\Lmin(\lambda_{1}\lambda_{2},\mu)$. Therefore
$\alpha\in\Ext(\lambda_{1}\lambda_{2};E)$, giving \[
\Ext(\lambda_{2};\Ext(\lambda_{1};E))\subseteq\Ext(\lambda_{1}\lambda_{2};E).\]

Now let $\alpha\in\Ext(\lambda_{1}\lambda_{2};E)$, so $(\alpha,\beta)\in\Lmin(\lambda_{1}\lambda_{2},\mu)$
for some $\mu\in E$ and $\beta\in\Lambda$. Then \begin{equation}
\lambda_{1}\lambda_{2}\alpha=\mu\beta~\text{ and }~d(\lambda_{1}\lambda_{2}\alpha)=d(\lambda_{1}\lambda_{2})\vee
d(\mu).\label{eqn:degree of alpha}\end{equation}
 It follows that \begin{align*}
\lambda_{2}\alpha & =(\mu\beta)(d(\lambda_{1}),d(\lambda_{1}\lambda_{2})\vee d(\mu))\\
 & =(\mu\beta)(d(\lambda_{1}),d(\lambda_{1})\vee d(\mu))(\mu\beta)(d(\lambda_{1})\vee
d(\mu),d(\lambda_{1}\lambda_{2})\vee d(\mu)).\end{align*}
 Defining $\xi:=(\mu\beta)(d(\lambda_{1}),d(\lambda_{1})\vee d(\mu))$,
we then have \[
\lambda_{2}\alpha=\xi(\mu\beta)(d(\lambda_{1})\vee d(\mu),d(\lambda_{1}\lambda_{2})\vee d(\mu)).\]
 To show that $\alpha\in\Ext(\lambda_{2};\Ext(\lambda_{1};E))$, we
need to show that $\xi\in\Ext(\lambda_{1};E)$ and $d(\lambda_{2}\alpha)=d(\lambda_{2})\vee d(\xi)$.
To this end, we calculate \begin{align}
\lambda_{1}\xi & =\lambda_{1}(\mu\beta)(d(\lambda_{1}),d(\lambda_{1})\vee d(\mu))\notag\\
 & =(\mu\beta)(0,d(\lambda_{1})\vee d(\mu))\quad\text{by \eqref{eqn:degree of alpha}}\label{eqn:degree of xi}\\
 & =\mu\beta(0,(d(\lambda_{1})\vee d(\mu))-d(\mu)).\notag\end{align}
 Furthermore, by \eqref{eqn:degree of xi} we have \begin{equation}
d(\lambda_{1}\xi)=d(\lambda_{1})\vee d(\mu).\label{eqn:degree of lambda1 xi}\end{equation}
 Hence $\xi\in\Ext(\lambda_{1};E)$. It remains to show that $d(\lambda_{2}\alpha)=d(\lambda_{2})\vee d(\xi)$.
On the one hand \begin{align*}
d(\lambda_{2}\alpha) & =(d(\lambda_{1}\lambda_{2})\vee d(\mu))-d(\lambda_{1})\quad\text{by \eqref{eqn:degree of
alpha}}\\
 & =d(\lambda_{2})\vee(d(\mu)-d(\lambda_{1})),\end{align*}
 and on the other hand \begin{align*}
d(\lambda_{2})\vee d(\xi) & =d(\lambda_{2})\vee((d(\lambda_{1})\vee d(\mu))-d(\lambda_{1}))\quad\text{by
\eqref{eqn:degree of lambda1 xi}}\\
 & =d(\lambda_{2})\vee(0\vee(d(\mu)-d(\lambda_{1})))\\
 & =d(\lambda_{2})\vee(d(\mu)-d(\lambda_{1})).\end{align*}
 So $d(\lambda_{2}\alpha)=d(\lambda_{2})\vee d(\xi)$, as required.
Therefore $\alpha\in\Ext(\lambda_{2};\Ext(\lambda_{1};E))$, and we
have \[
\Ext(\lambda_{2};\Ext(\lambda_{1};E))=\Ext(\lambda_{1}\lambda_{2};E).\]

\end{proof}

\section{Inverse Semigroups of Higher-Rank Graphs}

\label{sec:inverse semigroup}

For the remainder of the paper $(\Lambda,d)$ will be a finitely aligned
$k$-graph. To build the groupoids associated to $\Lambda$, we first construct an inverse semigroup $S_{\Lambda}$
associated to $\Lambda$ that captures the salient features of Toeplitz-Cuntz-Krieger
families. Recall that a semigroup $S$ is an \emph{inverse semigroup}
if for all $s\in S$ there exists a \emph{unique} element $s^{*}\in S$
such that $ss^{*}s=s$ and $s^{*}ss^{*}=s^{*}$. We follow \cite{P1} for the general theory and notation concerning
inverse semigroups. In particular, we denote the semilattice
of idempotents of $S$ by $E(S)$.

For $(\lambda,\mu),(\xi,\eta)\in\Lambda*_{s}\Lambda$,
we write $(\lambda,\mu)\perp(\xi,\eta)$ in case $\Lmin(\lambda,\xi)=\emptyset$
and $\Lmin(\mu,\eta)=\emptyset$.

\begin{dfn}\label{dfn:S_Lambda}
Define \(S_\Lambda\) to be the collection of all finite subsets \(F\) of \(\Lambda *_s \Lambda\) such that for
\emph{distinct} (\(\lambda\),\(\mu\)) and \((\nu,\omega)\in F\), we have \((\lambda,\mu)\perp (\nu,\omega)\).
\end{dfn}

\begin{rmk}
The empty subset of $\Lambda*_s\Lambda$ is an element of $S_\Lambda$.
\end{rmk}

\begin{prop}
For elements $F,G\in S_{\Lambda}$, the equation \begin{equation}
FG:=\bigcup_{(\lambda,\mu)\in F,(\xi,\eta)\in
G}\{(\lambda\alpha,\eta\beta):(\alpha,\beta)\in\Lmin(\mu,\xi)\}\label{eqn:def of S}\end{equation}
 defines an associative multiplication on $S_{\Lambda}$.
\end{prop}
\begin{proof}
The product of two finite sets is a finite set since $\Lambda$ is finitely aligned.
It is clear that the product also satisfies Definition~\ref{dfn:S_Lambda}, so $S_{\Lambda}$
is closed under multiplication.

For associativity, fix $(\lambda,\mu),(\xi,\eta),(\tau,\omega)\in\Lambda*_{s}\Lambda$.
It suffices to show that
\[
\big(\{(\lambda,\mu)\}\cdot\{(\xi,\eta)\}\big)\cdot\{(\tau,\omega)\}=\{(\lambda,\mu)\}\cdot
\big(\{(\xi,\eta)\}\cdot\{(\tau,\omega)\}\big).
\]
 The left-hand side equals
\[
\bigcup_{(\alpha,\beta)\in\Lmin(\mu,\xi)}\bigcup_{(\nu,\zeta)\in
\Lmin(\eta\beta,\tau)}\{(\lambda\alpha\nu,\omega\zeta)\}
\]
and the right-hand side equals
\[
\bigcup_{(\nu',\zeta')\in\Lmin(\eta,\tau)}\bigcup_{(\alpha',\beta')\in
\Lmin(\mu,\xi\nu')}\{(\lambda\alpha',\omega\zeta'\beta')\}.
\]
So we need to show that \begin{equation}
\begin{split}\{(\alpha\nu,\zeta):( & \alpha,\beta)\in\Lmin(\mu,\xi),(\nu,\zeta)\in\Lmin(\eta\beta,\tau)\}\\
 & =\{(\alpha',\zeta'\beta'):(\nu',\zeta')\in\Lmin(\eta,\tau),(\alpha',\beta')\in\Lmin(\mu,\xi\nu')\}.\end{split}
\label{eqn:associativity eqn}\end{equation}

Suppose $(\alpha\nu,\zeta)$ is an element of the left-hand side of
\eqref{eqn:associativity eqn}, so $(\alpha,\beta)\in\Lmin(\mu,\xi)$
and $(\nu,\zeta)\in\Lmin(\eta\beta,\tau)$. We will show that $(\alpha\nu,\zeta)$
is of the form $(\alpha',\zeta'\beta')$ where $(\nu',\zeta')\in\Lmin(\eta,\tau)$
and $(\alpha',\beta')\in\Lmin(\mu,\xi\nu')$. To begin, \begin{equation}
\mu\alpha=\xi\beta~\text{ and }~\eta\beta\nu=\tau\zeta\label{eqn:extension 1}\end{equation}
 where \begin{equation}
d(\mu\alpha)=d(\mu)\vee d(\xi)~\text{ and }~d(\eta\beta\nu)=d(\eta\beta)\vee d(\tau),\label{eqn:extension 1
degree}\end{equation}
 so we set \[
\nu':=(\beta\nu)(0,(d(\eta)\vee d(\tau))-d(\eta))~\text{ and }~\zeta':=\zeta(0,(d(\eta)\vee d(\tau))-d(\tau)).\]
 Then $(\nu',\zeta')\in\Lmin(\eta,\tau)$. By \eqref{eqn:extension 1},
\[
\mu\alpha\nu=\xi\beta\nu=\xi\nu'(\beta\nu)((d(\eta)\vee d(\tau))-d(\eta),d(\beta\nu)),\]
 so we set \[
\alpha':=\alpha\nu~\text{ and }~\beta':=(\beta\nu)((d(\eta)\vee d(\tau))-d(\eta),d(\beta\nu)).\]
 This gives $\mu\alpha'=\xi\nu'\beta'$ and $\zeta'\beta'=\zeta$.

We now have $(\alpha\nu,\zeta)=(\alpha',\zeta'\beta')$ where $(\nu',\zeta')\in\Lmin(\eta,\tau)$
and $\mu\alpha'=\xi\nu'\beta'$, so it remains to show that $(\alpha',\beta')\in\Lmin(\mu,\xi\nu')$;
that is, $d(\mu\alpha')=d(\mu)\vee d(\xi\nu')$. For this we calculate
\begin{align*}
d(\mu\alpha') & =d(\mu\alpha)+d(\nu)\\
 & =\big(d(\mu)\vee d(\xi)\big)+d(\nu) \quad\text{by \eqref{eqn:extension 1 degree}}\\
 & =\big(d(\mu)\vee d(\xi)\big)+\big(\big(d(\eta\beta)\vee d(\tau)\big)-d(\eta\beta)\big)\quad\text{by
\eqref{eqn:extension 1 degree}}\\
 & =\big(d(\mu)\vee d(\xi)\big)+\big(0\vee\big(d(\tau)-d(\eta)-d(\beta)\big)\big)\\
 & =\big(d(\mu)\vee d(\xi)\big)+\big(0\vee\big(d(\tau)-d(\eta)-\big(d(\mu)\vee
d(\xi)\big)+d(\xi)\big)\big)\quad\text{by \eqref{eqn:extension 1 degree}}\\
 & =d(\mu)\vee d(\xi)\vee\big(d(\tau)-d(\eta)+d(\xi)\big)\\
 & =d(\mu)\vee\big(d(\xi)+\big(d(\eta)\vee d(\tau)\big)-d(\eta)\big)\\
 & =d(\mu)\vee\big(d(\xi)+d(\nu')\big)\quad\text{since }(\nu',\zeta')\in\Lmin(\eta,\tau)\\
 & =d(\mu)\vee d(\xi\nu'),\end{align*}
 as required. Therefore $(\alpha\nu,\zeta)$ is of the form $(\alpha',\zeta'\beta')$
where $(\nu',\zeta')\in\Lmin(\eta,\tau)$ and $(\alpha',\beta')\in\Lmin(\mu,\xi\nu')$,
and so the left-hand side of \eqref{eqn:associativity eqn} is contained
in the right-hand side of \eqref{eqn:associativity eqn}.

The reverse containment can be proved in a similar way, giving the
result.
\end{proof}
\begin{prop}
$S_{\Lambda}$ is an inverse semigroup with involution defined for
$F\in S_{\Lambda}$ by $F^{*}:=\{(\mu,\lambda):(\lambda,\mu)\in F\}$.
\end{prop}
\begin{proof}
Fix $F\in S_{\Lambda}$. We need to show that $F^{*}$ is the unique
element of $S_{\Lambda}$ such that $FF^{*}F=F$ and $F^{*}FF^{*}=F^{*}$.
Using \eqref{eqn:def of S} we have \[
F(F^{*}F)=F\{(\mu,\mu):(\lambda,\mu)\in F\}=F,\]
 and similarly, $F^{*}FF^{*}=F^{*}$. For uniqueness, suppose $FGF=F$
and $GFG=G$ for some $G\in S_{\Lambda}$. The product $FGF$ can
be written as \[
\bigcup_{(\lambda,\mu),(\tau,\omega)\in F}\bigcup_{(\xi,\eta)\in
G}\{(\lambda\alpha\nu,\omega\zeta):(\alpha,\beta)\in\Lmin(\mu,\xi),(\nu,\zeta)\in\Lmin(\eta\beta,\tau)\}.\] Since
$FGF=F$, for each $(\lambda',\mu')\in F$, there exist $(\lambda,\mu),(\tau,\omega)\in F$,
$(\xi,\eta)\in G$, $(\alpha,\beta)\in\Lmin(\mu,\xi)$ and $(\nu,\zeta)\in\Lmin(\eta\beta,\tau)$
such that $(\lambda',\mu')=(\lambda\alpha\nu,\omega\zeta)$. This shows that $(\lambda',\mu')\not\perp(\lambda,\mu)$ and
$(\lambda',\mu')\not\perp(\tau,\omega)$.  Thus Definition \ref{dfn:S_Lambda} implies
$(\lambda',\mu')=(\lambda,\mu)=(\tau,\omega)$,
so $\alpha=\nu=\zeta=s(\lambda)$. Hence $(s(\lambda),\beta)\in\Lmin(\mu,\xi)$
and $(s(\lambda),s(\lambda))\in\Lmin(\eta\beta,\lambda)$, which give
\begin{equation}
\lambda=\eta\beta~\text{ and }~\mu=\xi\beta.\label{eqn:mu and lambda extensions}\end{equation}

We claim that $\beta=s(\xi)$. Since $(\xi,\eta)\in G$, $(\lambda,\mu)\in F$,
$(\beta,s(\lambda))\in\Lmin(\eta,\lambda)$ and $(s(\lambda),\beta)\in\Lmin(\mu,\xi)$,
it follows from the definition of multiplication that $(\xi\beta,\eta\beta)\in GFG$.
Since $GFG=G$, it then follows from Definition~\ref{dfn:S_Lambda} that $(\xi\beta,\eta\beta)=(\xi,\eta)$,
which gives $\beta=s(\xi)$. Therefore \eqref{eqn:mu and lambda extensions}
gives $\lambda=\eta$ and $\mu=\xi$. Thus for every $(\lambda,\mu)\in F$,
we have $(\mu,\lambda)\in G$, so $F^{*}\subseteq G$. A similar argument
gives $G^{*}\subseteq F$ which is equivalent to $G\subseteq F^{*}$.
Therefore $G=F^{*}$, and $S_{\Lambda}$ is an inverse semigroup.
\end{proof}
\begin{rmk}
As we noted in Section~\ref{sec:Motivation}, the inverse semigroup Paterson used in \cite{P2} is $\Lambda*_s \Lambda \cup \{z\}$.  The multiplication defined on Paterson's inverse semigroup corresponds to the multiplication of singleton sets
in $S_\Lambda$. If $\Lambda$ is a $k$-graph with $k\geq 2$ and $\lambda,\mu\in\Lambda$, then $\Lmin(\lambda,\mu)$
generally contains more than one element.  Thus $S_\Lambda$ was chosen to consist of finite subsets of
$\Lambda*_s\Lambda$ in order to obtain an inverse semigroup that is closed under a multiplication and which captures
the properties of the $k$-graph.  We note, too, that in \cite{P2} one has to add a zero element to
$\Lambda*_s \Lambda$, but this is not necessary in our setting, since the empty subset of $\Lambda*_s \Lambda$ is
included in our $S_{\Lambda}$ and serves as a zero.
\end{rmk}

\section{The action of $S_{\Lambda}$ on the path space $X_{\Lambda}$}

\label{sec:the action}

\begin{dfn}
Let \(X_\Lambda\) be the set of all graph morphisms \(x:\Omega_{k,m}\to\Lambda\) where \(m\in(\NN\cup\{\infty\})^k\)
(See Definition~\ref{dfn:morphism} and Example~\ref{k-graph examples}). We extend range and degree maps to all
\(x:\Omega_{k,m}\to\Lambda\) in \(X_\Lambda\) by setting \(r(x):=x(0)\) and \(d(x):=m\).
\end{dfn}

\begin{notation}
\noindent Let $x\in X_{\Lambda}$.
\begin{enumerate}
\item  For $m\in\NN^{k}$ with $m\leq d(x)$, we
write $\bp(x;m,{d(x)})$ for the graph morphism $\bp(x;m,{d(x)}):\Omega_{k,d(x)-m}\rightarrow\Lambda$
defined by the equation \[
\bp(x;m,{d(x)})(p,q)=x(p+m,q+m)\text{ for all $p,q\in\NN^{k}$, $p\leq q\leq d(x)-m$.}\]
\item For $\lambda\in\Lambda x(0)$, we write $\lambda x$ for the graph
morphism $\lambda x:\Omega_{k,d(x)+d(\lambda)}\rightarrow\Lambda$
defined by $(\lambda x)(0,m)=\lambda x(0,m-d(\lambda))$ for all $m\in\NN^{k}$,
$d(\lambda)\leq m\leq d(x)+d(\lambda)$.
\end{enumerate}
\end{notation}

\begin{rmk}
The set $X_\Lambda$ contains $\Lambda$ as a subset if we identify each path $\lambda\in\Lambda$ with the unique graph
morphism $x_\lambda:\Omega_{k,d(\lambda)}\to\Lambda$ given by $x_\lambda(0,d(\lambda))=\lambda$.
\end{rmk}

For $F\in S_{\Lambda}$ define \[
D_{F}:=\{ x\in X_{\Lambda}:\text{ there exists }(\lambda,\mu)\in F\text{ such that }x(0,d(\mu))=\mu\}\]
 and define $\theta_{F}:D_{F}\rightarrow D_{F^{*}}$ by $\theta_{F}(x):=\lambda\bp(x;{d(\mu)},{d(x)})$
where $(\lambda,\mu)$ is the unique element of $F$ such that $x(0,d(\mu))=\mu$.
For $\lambda\in\Lambda$, we write $D_{\lambda}$ instead of $D_{\{(\lambda,\lambda)\}}$.

\begin{prop}
\label{prop:basis for X} The family of sets $\{ D_{F},D_{F}^{c}:F\in S_{\Lambda}\}$
is a subbasis for a Hausdorff topology on $X_{\Lambda}$.
\end{prop}
\begin{proof}
Since $X_{\Lambda}=\bigcup_{v\in\Lambda^{0}}D_{v}$, the family $\{ D_{F},D_{F}^{c}:F\in S_{\Lambda}\}$
forms a subbasis for a topology on $X_{\Lambda}$. To show that the
topology is Hausdorff, fix $x,y\in X_{\Lambda}$, $x\neq y$. If $x(0)\neq y(0)$,
then $x\in D_{x(0)}$, $y\in D_{y(0)}$ and $D_{x(0)}\cap D_{y(0)}=\emptyset$.
So assume that $x(0)=y(0)$, and let $n\in\NN^{k}$ be minimal with
respect to the condition \[
x(0,n)=y(0,n)\text{ and }x(0,n+e_{i})\neq y(0,n+e_{i})\text{ for some }i\in\{1,\dots,k\}.\]
 Without loss of generality assume that $n+e_{i}\leq d(x)$. Then
$x\in D_{x(0,n+e_{i})}$ and $y\in D_{x(0,n+e_{i})}^{c}$, as required.
\end{proof}
\begin{rmk}\label{rmk:basis for topology on X}
For \(F_1,\dots,F_n,G_1,\dots,G_m\in S_\Lambda\), letting \(F=\prod_{i=1}^nF_i^*F_i\), we have
\[
D_{F_1}\cap\cdots\cap D_{F_n}\cap D_{G_1}^c\cap\cdots\cap D^c_{G_m} =
D_F\cap\textstyle{\bigcap^m_{j=1}D^c_{FG_j^*G_j}}.
\]
Since \(D_F\) is the disjoint union of the \(D_\lambda\) where \((\lambda,\lambda)\in F\), the family of sets
\[
\{D_\lambda\cap D^c_{\lambda\nu_1}\cap\cdots\cap D^c_{\lambda\nu_l}: \lambda\in\Lambda \text{ and }
\nu_1,\dots,\nu_l\in s(\lambda)\Lambda\}
\]
is a basis for the given topology on \(X_\Lambda\).
\end{rmk}

\begin{rmk}
An infinite sequence of paths in \(\Lambda\) is called \emph{wandering} if for any finite set \(E\subseteq\Lambda\),
the sequence is eventually in \(\Lambda\setminus E\). Convergence in \(X_\Lambda\) is then given by: A sequence
$\langle x_i\rangle$ converges to $x$ if and only if the following two conditions occur:
\begin{itemize}
\item[(1)]
for all \(n\in\NN^k\) such that \(n\le d(x)\), there exists \(I\in\NN\) such that \(i\ge I\) implies
\(x_i(0,n)=x(0,n)\), and
\item[(2)]
if \(d(x)_j<\infty\), then for any \(n\in\NN^k\) such that \(n\le d(x)\) and \(n_j=d(x)_j\), either
\[
J(j,n):=\{i\in\NN : d(x_i)\ge n+e_j\}
\]
is finite, or \(\langle x_i(n,n+e_j)\rangle_{i\in J(j,n)}\) is wandering.
\end{itemize}
\end{rmk}

\begin{prop}
\label{prop:D_v compact} For each $v\in\Lambda^{0}$, $D_{v}$ is
compact, and so, in particular, $X_{\Lambda}$ is a locally compact Hausdorff space.
\end{prop}
\begin{proof}
Let $\langle x_{i}\rangle_{i=1}^{\infty}$ be a sequence in $D_{v}$.
We construct an element $x\in X_{\Lambda}$ and a subsequence $\langle y_{i}\rangle_{i=1}^{\infty}$
such that $y_{i}\rightarrow x$.

Let $\{ E_{n}\}_{n=1}^{\infty}$ be a listing of all finite subsets
of $v\Lambda$ containing $v$. There is at least one $\lambda\in E_{1}$
such that $x_{i}(0,d(\lambda))=\lambda$ for infinitely many $i\in\NN$
(namely $\lambda=v$). Let $\lambda_{1}$ be such a $\lambda$ of
maximal degree.

Suppose that $\lambda_{1},\dots,\lambda_{n}$ have already been defined
so that $\lambda_{1}\cdots\lambda_{n}\in\,\, v\Lambda$ and
$x_{i}(0,d(\lambda_{1}\cdots\lambda_{n}))=\lambda_{1}\cdots\lambda_{n}$
for infinitely many $i\in\NN$. There is at least one $\lambda\in\Ext(\lambda_{1}\cdots\lambda_{n};E_{n+1})$
such that $x_{i}(0,d(\lambda_{1}\cdots\lambda_{n}\lambda))=\lambda_{1}\cdots\lambda_{n}\lambda$
for infinitely many $i\in\NN$ (namely $\lambda=s(\lambda_{n})$,
since $v\in E_{n+1}$). Let $\lambda_{n+1}$ be such a $\lambda$
of maximal degree.

Define $m\in(\NN\cup\{\infty\})^{k}$ by $m:=\lim_{n\rightarrow\infty}d(\lambda_{1}\cdots\lambda_{n})$.
There is a unique $x\in X_{\Lambda}$ such that $d(x)=m$ and
$x(0,d(\lambda_{1}\cdots\lambda_{n}))=\lambda_{1}\cdots\lambda_{n}$
for all $n\in\NN$. Let $\langle y_{i}\rangle_{i=1}^{\infty}$ be
a subsequence of $\langle x_{i}\rangle_{i=1}^{\infty}$ such that
for all $n\in\NN$, \begin{equation}
y_{i}(0,d(\lambda_{1}\cdots\lambda_{n}))=\lambda_{1}\cdots\lambda_{n}\text{ for all }i\geq n.\label{eqn:def of
y_i}\end{equation}
 We claim that $\lim_{i\rightarrow\infty}y_{i}=x$.

Fix a neighborhood $D_{\lambda}\cap D_{\lambda\nu_{1}}^{c}\cap\cdots\cap D_{\lambda\nu_{l}}^{c}$
of $x$, so we have $x(0,d(\lambda))=\lambda$ and $x(0,d(\lambda\nu_{j}))\neq\lambda\nu_{j}$
for $j=1,\dots,l$. There exists $n\in\NN$ such that \[
E_{n}=\{ v\}\cup\{\lambda\nu_{1},\dots,\lambda\nu_{l}\}\]
 and \begin{align*}
\lambda_{n}\in & \Ext(\lambda_{1}\cdots\lambda_{n-1};E_{n})\\
& =\{ s(\lambda_{n-1})\}\cup\textstyle{\bigcup_{j=1}^{l}\{\alpha:(\alpha,\beta)\in
\Lmin(\lambda_{1}\cdots\lambda_{n-1},\lambda\nu_{j}) \text{ for some } \beta\in \Lambda\}}.\end{align*}

For $i\geq n-1$, if $y_{i}(0,d(\lambda\nu_{j}))=\lambda\nu_{j}$
for some $j\in\{1,\dots,l\}$, then by \eqref{eqn:def of y_i}\begin{equation}
\begin{split}y_{i}(d(\lambda_{1}\cdots\lambda_{n-1}), & d(\lambda_{1}\cdots\lambda_{n-1})\vee d(\lambda\nu_{j}))\\
\in & \{\alpha:(\alpha,\beta)\in\Lmin(\lambda_{1}\cdots\lambda_{n-1},\lambda\nu_{j}) \text{ for some }
\beta\in\Lambda\}.\end{split}
\label{eqn:sequence in set}\end{equation}

Suppose, for contradiction, that there are infinitely many $i\geq n-1$
such that $y_{i}(0,d(\lambda\nu_{j}))=\lambda\nu_{j}$ for some $j\in\{1,\dots,l\}$.
Then there are infinitely many $i\geq n-1$ such that \eqref{eqn:sequence in set}
holds for some $j\in\{1,\dots,l\}$. Since $\lambda_{n}$ is chosen
to be of maximal degree, we must then have \[
\lambda_{n}\in\bigcup_{j=1}^{l}\{\alpha:(\alpha,\beta)\in\Lmin(\lambda_{1}\cdots\lambda_{n-1},\lambda\nu_{j})\}\]
 But then $x(0,d(\lambda_{1}\cdots\lambda_{n}))=\lambda_{1}\cdots\lambda_{n}=\lambda\nu_{j}\beta$
for some $j\in\{1,\dots,l\}$ and $\beta\in\Lambda$, contradicting
the inequality $x(0,d(\lambda\nu_{j}))\neq\lambda\nu_{j}$. Hence
there are only finitely many $i\in\NN$ such that $y_{i}\not\in D_{\lambda}\cap D_{\lambda\nu_{1}}^{c}\cap\cdots\cap
D_{\lambda\nu_{l}}^{c}$.
Therefore $\lim_{i\rightarrow\infty}y_{i}=x$, and $D_{v}$ is compact.
\end{proof}
\begin{cor}
\label{cor:D_F compact} For each $F\in S_{\Lambda}$, $D_{F}$ is
compact.
\end{cor}
\begin{proof}
For each $\lambda\in\Lambda$, $D_{\lambda}$ is a closed subset of
$D_{r(\lambda)}$, and hence is compact by Proposition~\ref{prop:D_v
compact}. The result then follows since $D_{F}=\bigcup_{(\lambda,\mu)\in F}D_{\mu}$
-- a finite union.
\end{proof}
\begin{rmk}\label{rmk:two groupoids rmk}
In the next section we will construct two topological groupoids associated to a finitely aligned \(k\)-graph
\(\Lambda\). The unit space of the first groupoid \(\Gg_\Lambda\) is homeomorphic to \(X_\Lambda\), and
\(C^*(\Gg_\Lambda)\) is isomorphic to the Toeplitz algebra \(\Tt C^*(\Lambda)\). The second groupoid is a reduction of
\(\Gg_\Lambda\) to a closed invariant subset \(\partial\Lambda\) of \((\Gg_\Lambda)^{(0)}\); we describe
\(\partial\Lambda\) as a subspace of \(X_\Lambda\). The \(C^*\)-algebra of \(\Gg_\Lambda|_{\partial\Lambda}\) is
isomorphic to the Cuntz-Krieger algebra \(C^*(\Lambda)\).
\end{rmk}

\begin{dfn}
An element \(x\in X_\Lambda\) is called a \emph{boundary path} if for all \(n\in\NN^k\) with \(n\le d(x)\), and for all
\(E\in x(n)\FE(\Lambda)\), there exists \(\lambda\in E\) such that \(x(n,n+d(\lambda))=\lambda\). We write
\(\partial\Lambda\) for the set of all boundary paths, and for \(v\in\Lambda^0\), write \(v(\partial\Lambda)\) for
\(\{x\in\partial\Lambda  :  r(x)=v\}\).
\end{dfn}

\begin{examples}\label{ex:boundary paths}
\begin{enumerate}
\item If $\Lambda$ is a row-finite $k$-graph with no sources, $\partial\Lambda=\Lambda^{\infty}$ where
$\Lambda^{\infty}$ is the set of all graph morphisms from $\Omega_{k,(\infty,\infty,\ldots,\infty)}$ to $\Lambda$.
\item For the $k$-graph in Example \ref{ex:finitely aligned} (2), any path starting with $\gamma,\eta,\xi,\omega$ or
$\tau_i$ for some $i$ is a boundary path. Furthermore, since any finite exhaustive subset of $w\Lambda$ must contain
$w$, every element in $\Lambda w$ is also a boundary path.
\end{enumerate}
\end{examples}

\begin{lemma}\label{lem:partialLambda closed}
\(\partial\Lambda\) is closed in \(X_\Lambda\).
\end{lemma}

\begin{proof}
$\partial\Lambda$ is closed since it is the complement in $X_\Lambda$ of the open set
\[
\bigcup_{\lambda\in\Lambda} \bigcup_{E\in s(\lambda)\FE(\Lambda)} D_\lambda \cap \textstyle{\bigcap_{\nu\in E}
D_{\lambda\nu}^c}.
\]
\end{proof}

Straightforward calculations give the following lemma.
\begin{lemma}\label{lem:partialLambda closed under action}
Let \((\Lambda,d)\) be a finitely aligned \(k\)-graph and let \(x\in\partial\Lambda\).
\begin{enumerate}
\item
If \(m\in\NN^k\) and \(m\le d(x)\), then \(\bp(x;m,{d(x)})\in\partial\Lambda\).
\item
If \(\lambda\in\Lambda x(0)\), then \(\lambda x\in \partial\Lambda\).
\end{enumerate}
\end{lemma}

\begin{rmk}\label{rmk:partial Lambda closed invariant}
Lemma~\ref{lem:partialLambda closed} and Lemma~\ref{lem:partialLambda closed under action} imply that
\(\partial\Lambda\) is a locally compact Hausdorff space which is invariant under the action of \(S_\Lambda\).
Following on from Remark~\ref{rmk:two groupoids rmk}, this fact will give us the required closed invariant subset of
the unit space of \(\Gg_\Lambda\). The next lemma ensures that the partial isometries inside
\(C^*(\Gg_\Lambda|_{\partial\Lambda})\) giving a Cuntz-Krieger \(\Lambda\)-family are nonzero.
\end{rmk}

\begin{lemma}\label{lem:v partialLambda nonempty}
For all \(v\in\Lambda^0\), \(v(\partial\Lambda)\) is nonempty.
\end{lemma}

\begin{proof}
We construct a sequence of paths $\lambda_{1},\lambda_{2},\dots$
such that $r(\lambda_{1})=v$ and $s(\lambda_{n})=r(\lambda_{n+1})$.
We then show that the sequence of paths defines an element of $v(\partial\Lambda)$.
\begin{enumerate}
\item Let $\{ E_{1,i}\}_{i=1}^{\infty}$ be a listing of $v\FE(\Lambda)$.
Choose $\lambda_{1}\in E_{1,1}$.
\item Let $\{ E_{2,i}\}_{i=1}^{\infty}$ be a listing of $s(\lambda_{1})\FE(\Lambda)$.
Choose $\lambda_{2}\in\Ext(\lambda_{1};E_{1,2})$.
\item Let $\{ E_{3,i}\}_{i=1}^{\infty}$ be a listing of $s(\lambda_{2})\FE(\Lambda)$.
Choose $\lambda_{3}\in\Ext(\lambda_{2};E_{2,1})$.
\item Let $\{ E_{4,i}\}_{i=1}^{\infty}$ be a listing of $s(\lambda_{3})\FE(\Lambda)$.
Choose $\lambda_{4}\in\Ext(\lambda_{1}\lambda_{2}\lambda_{3};E_{1,3})$.
\item Let $\{ E_{5,i}\}_{i=1}^{\infty}$ be a listing of $s(\lambda_{4})\FE(\Lambda)$.
Choose $\lambda_{5}\in\Ext(\lambda_{2}\lambda_{3}\lambda_{4};E_{2,2})$.
\item Let $\{ E_{6,i}\}_{i=1}^{\infty}$ be a listing of $s(\lambda_{5})\FE(\Lambda)$.
Choose $\lambda_{6}\in\Ext(\lambda_{3}\lambda_{4}\lambda_{5};E_{3,1})$.
\end{enumerate}
Continuing this process, we define $\lambda_{n}\in\Lambda$ for all
$n\in\NN$ satisfying $r(\lambda_{1})=v$, $s(\lambda_{n})=r(\lambda_{n+1})$
and for all $E\in s(\lambda_{n})\FE(\Lambda)$ there exists $m\in\NN$
such that $\lambda_{m}\in\Ext(\lambda_{n+1}\cdots\lambda_{m-1};E)$.
There exists a unique graph morphism $x\in X_{\Lambda}$ such that
$d(x)=\lim_{n\rightarrow\infty}d(\lambda_{1}\cdots\lambda_{n})$ and
$x(0,d(\lambda_{1}\cdots\lambda_{n}))=\lambda_{1}\cdots\lambda_{n}$
for all $n\in\NN$.

Let $p\in\NN^{k}$, $p\leq d(x)$, and let $E\in x(p)\FE(\Lambda)$.
We will show that there exists $\mu\in E$ such that $x(p,p+d(\mu))=\mu$.
Let $n\in\NN$ be the smallest number such that $d(\lambda_{1}\cdots\lambda_{n})\geq p$.
Then $\Ext(x(p,d(\lambda_{1}\cdots\lambda_{n}));E)\in s(\lambda_{n})\FE(\Lambda)$
by Proposition~\ref{prop:C.4 RSY2}, so there exists $i\in\NN$ such
that \[
E_{n+1,i}=\Ext(x(p,d(\lambda_{1}\cdots\lambda_{n}));E).\]
 By construction of $x$, there exists $m\in\NN$ such that \begin{align*}
\lambda_{m}\in & \Ext(\lambda_{n+1}\cdots\lambda_{m-1};E_{n+1,i})\\
 & =\Ext\big(\lambda_{n+1}\cdots\lambda_{m-1};\Ext(x(p,d(\lambda_{1}\cdots\lambda_{n}));E)\big)\\
 & =\Ext(x(p,d(\lambda_{1}\cdots\lambda_{m-1}));E)\quad\text{by Proposition~\ref{prop:Ext transitive}.}\end{align*}
 Thus \begin{align*}
x(p,d(\lambda_{1}\cdots\lambda_{m})) & =x(p,d(\lambda_{1}\cdots\lambda_{m-1}))\lambda_{m}\\
 & =\mu\alpha\end{align*}
 for some $\mu\in E$ and $\alpha\in\Lambda$, so $x(p,p+d(\mu))=\mu$.
Therefore $x\in v(\partial\Lambda)$, as required.
\end{proof}

\section{The groupoid of the system $(X_{\Lambda},S_{\Lambda},\theta)$}

\label{sec:the groupoid}

We begin this section by recalling notions from groupoid theory. We
follow Paterson's book \cite{P1}. Let $G$ be a topological groupoid.
The unit space of $G$ is denoted $G^{(0)}$ and the space of composable
pairs is denoted $G^{(2)}$. We use $r$ and $s$ for the range and
source maps $r(g)=gg^{-1}$ and $s(g)=g^{-1}g$ (for $g\in G$). A
subset $U\subseteq G^{(0)}$ is invariant if $r(s^{-1}(U))$ is contained
in $U$. We denote by $G^{op}$ the family of open Hausdorff subsets
$A$ of $G$ such that $r|_{A}$ and $s|_{A}$ are homeomorphisms.
A locally compact groupoid $G$ is called \emph{$r$-discrete} if
$G^{(0)}$ is open in $G$. In fact, a locally compact groupoid $G$
is $r$-discrete and admits a left Haar system if and only if $G^{op}$
is a basis for the topology of $G$, which happens if and only if
$r:G\rightarrow G^{(0)}$ is a local homeomorphism (see \cite[Proposition~2.8]{R}).
We define \[
G^{a}:=\{ A\in G^{op}:A\text{ is compact}\}.\]
 An $r$-discrete groupoid $G$ is called \emph{ample} if $G^{a}$
forms a basis for the topology of $G$.

We now construct the groupoid of germs of the system $(X_{\Lambda},S_{\Lambda},\theta)$.
Define \[
\Xi_{\Lambda}:=\{(F,x):F\in S_{\Lambda},x\in D_{F}\}\]
 and define a relation $\sim$ on $\Xi_{\Lambda}$ by requiring that
$(F,x)\sim(G,y)$ if and only if $x=y$ and there exists $P\in E(S_{\Lambda})$
such that $x\in D_{P}$ and $FP=GP$. \begin{lemma}
The relation \(\sim\) is an equivalence relation on \(\Xi_\Lambda\).
\end{lemma}

\begin{proof}
The reflexivity and symmetry of $\sim$ are obvious. Suppose $(F,x)\sim(G,x)$
and $(G,x)\sim(H,x)$. Then there exist $P_{1},P_{2}\in E(S_{\Lambda})$
such that $x\in D_{P_{1}}$, $x\in D_{P_{2}}$, $FP_{1}=GP_{1}$ and
$GP_{2}=HP_{2}$. Thus $x\in D_{P_{1}P_{2}}$ and \[
FP_{1}P_{2}=GP_{1}P_{2}=GP_{2}P_{1}=HP_{2}P_{1}=HP_{1}P_{2},\]
 so $(F,x)\sim(H,x)$, and $\sim$ is an equivalence relation.
\end{proof}

The equivalence class of $(F,x)$ will be denoted $[F,x]$. It is
called the \emph{germ} of $F$ at $x$.

Write $\Gg_{\Lambda}$ for the set $\Xi_{\Lambda}/\sim$. Then $\Gg_{\Lambda}$
becomes a groupoid where the composable pairs are of the form $\big([F,\theta_{G}(x)],[G,x]\big)$,
and product and inversion are given by the formulas \[
[F,\theta_{G}(x)]\cdot[G,x]=[FG,x]\quad\text{and}\quad[F,x]^{-1}=[F^{*},\theta_{F}(x)].\]
 The unit space $\Gg_{\Lambda}^{(0)}$ will be identified with $X_{\Lambda}$
via the map $[\{(r(x),r(x))\},x]\leftrightarrow x$. Then for $[F,x]\in\Gg_{\Lambda}$,
we have $r([F,x])=\theta_{F}(x)$ and $s([F,x])=x$.

\begin{rmk}\label{rmk:composable pairs}
For \(F\in S_\Lambda\) and \(x\in D_F\), let \((\lambda,\mu)\in F\) be the unique element of \(F\) such that
\(x(0,d(\mu))=\mu\). Then \([F,x]=[\{(\lambda,\mu)\},x]\), so \(\Gg_\Lambda\) can be expressed as
\[
\Gg_\Lambda = \{[\{(\lambda,\mu)\},x]  :  (\lambda,\mu)\in\Lambda*_s\Lambda, x\in D_\mu\}.
\]
Using this description of \(\Gg_\Lambda\), for \([\{(\lambda,\mu)\},x],[\{(\xi,\eta)\},x]\in \Gg_\Lambda\),
\begin{equation}\label{eqn:equivalence rel}
\begin{split}
[\{(\lambda,\mu)\},x]&=[\{(\xi,\eta)\},x] \\
\iff \lambda x(d(\mu),d(\mu)\vee d(&\eta))=\xi x(d(\eta),d(\mu)\vee d(\eta)).
\end{split}
\end{equation}
Furthermore, if \(\left([\{(\lambda,\mu)\}, x],[\{(\xi,\eta)\},y]\right)\in\Gg_\Lambda^{(2)}\), then setting
\[
\alpha:=x(d(\mu),d(\mu)\vee d(\xi)) \text{ and } \beta:=\bp(y;{d(\eta)},{d(y)})(0,(d(\mu)\vee d(\xi))-d(\xi)),
\]
we have \((\alpha,\beta)\in\Lmin(\mu,\xi)\) and
\begin{equation}\label{eqn:composable pairs}
[\{(\lambda,\mu)\}, x] \cdot [\{(\xi,\eta)\},y] = [\{(\lambda\alpha,\eta\beta)\},y]
\end{equation}
\end{rmk}

For $F\in S_{\Lambda}$, define \begin{equation}
\Psi(F):=\{[F,x]:x\in D_{F}\}\subseteq\Gg_{\Lambda}.\label{eqn:Psi}\end{equation}

\begin{prop}

\label{prop:G is an S-groupoid} The set $\{\Psi(F),\Psi(F)^{c}:F\in S_{\Lambda}\}$
is a subbasis for a locally compact Hausdorff topology on $\Gg_{\Lambda}$
in which each $\Psi(F)$ is compact.
\end{prop}
\begin{proof}
Since $\Gg_{\Lambda}=\bigcup_{(\lambda,\mu)\in\Lambda*_{s}\Lambda}\Psi(\{(\lambda,\mu)\})$,
the set $\{\Psi(F),\Psi(F)^{c}:F\in S_{\Lambda}\}$ is a subbasis
for a topology on $\Gg_{\Lambda}$.

For $(\lambda,\mu)\in\Lambda*_{s}\Lambda$, the source map in $\Gg_{\Lambda}$
restricted to $\Psi(\{(\lambda,\mu)\})$ is injective with image $D_{\mu}$.
The topology on $D_{\mu}$ is generated by a subbasis comprising elements
of the form $D_{\mu\alpha}$ and $D_{\mu}\cap D_{\mu\alpha}^{c}$.
The inverse images of these subbasis elements, \[
\left(s|_{\Psi(\{(\lambda,\mu)\})}\right)^{-1}(D_{\mu\alpha})=\Psi(\{(\lambda\alpha,\mu\alpha)\})\]
 and \[
\left(s|_{\Psi(\{(\lambda,\mu)\})}\right)^{-1}(D_{\mu}\cap
D_{\mu\alpha}^{c})=\Psi(\{(\lambda,\mu)\})\cap\Psi(\{(\lambda\alpha,\mu\alpha)\})^{c},\]
 form a subbasis for the topology on $\Psi(\{(\lambda,\mu)\})$. Therefore
$s|_{\Psi(\{(\lambda,\mu)\})}$ is a homeomorphism, and so Proposition~\ref{prop:basis for X}
and Corollary~\ref{cor:D_F compact} imply that $\Psi(\{(\lambda,\mu)\})$
is Hausdorff and compact. The result follows.
\end{proof}
\begin{rmk}
The family of sets of the form
\[
\Psi(\{(\lambda,\mu)\})\cap \Psi(\{(\lambda\nu_1,\mu\nu_1)\})^c\cap\cdots\cap \Psi(\{(\lambda\nu_l,\mu\nu_l)\})^c,
\]
where \((\lambda,\mu)\in\Lambda *_s\Lambda\) and \(\nu_1,\dots,\nu_l\in s(\lambda)\Lambda\), is a basis for the
topology on \(\Gg_\Lambda\).
\end{rmk}

\begin{prop}
$\Gg_{\Lambda}$ is an $r$-discrete topological groupoid.
\end{prop}
\begin{proof}
To show that composition is continuous, fix a composable pair \[
\left([\{(\lambda,\mu)\},x],[\{(\xi,\eta)\},y]\right)\in\Gg_{\Lambda}^{(2)}.\]
 Remark~\ref{rmk:composable pairs} says that the composable pair
is of the form \[
\left([\{(\lambda\alpha,\mu\alpha)\},\xi\beta z],[\{(\xi\beta,\eta\beta)\},\eta\beta z]\right),\]
 where $(\alpha,\beta)\in\Lmin(\mu,\xi)$, and has product $[\{(\lambda\alpha,\eta\beta)\},\eta\beta z]$.
Let
\[
A:=\Psi(\{(\lambda\alpha\tau,\eta\beta\tau)\})\cap\textstyle{\bigcap_{i=1}^{l}
\Psi(\{(\lambda\alpha\tau\nu_{i},\eta\beta\tau\nu_{i})\})^{c}}
\]
be a neighborhood of $[\{(\lambda\alpha,\eta\beta)\},\eta\beta z]$.
Then for neighborhoods
\[
B:=\Psi(\{(\lambda\alpha\tau,\mu\alpha\tau)\})\cap
\textstyle{\bigcap_{i=1}^{l}\Psi(\{(\lambda\alpha\tau\nu_{i},\mu\alpha\tau\nu_{i})\})^{c}}
\]
and
\[
C:=\Psi(\{(\xi\beta\tau,\eta\beta\tau)\})\cap
\textstyle{\bigcap_{i=1}^{l}\Psi(\{(\xi\beta\tau\nu_{i},\eta\beta\tau\nu_{i})\})^{c}}
\]
of $[\{(\lambda\alpha,\mu\alpha)\},\xi\beta z]$ and $[\{(\xi\beta,\eta\beta)\},\eta\beta z]$,
respectively, we have $BC\subseteq A$, which gives continuity of
composition.

For $F\in S_{\Lambda}$, $\Psi(F)=\Psi(F^{*})^{-1}$ and $\Psi(F)^{c}=\left(\Psi(F^{*})^{c}\right)^{-1}$,
so inversion is continuous.

Finally, $\Gg_{\Lambda}$ is $r$-discrete since $\Gg_{\Lambda}^{(0)}=\bigcup_{v\in\Lambda^{0}}\Psi(\{(v,v)\})$
is open in $\Gg_{\Lambda}$.
\end{proof}
\begin{rmk}\label{rmk:unit space topology}
The identification of \(X_\Lambda\) with \(\Gg_\Lambda^{(0)}\) is a homeomorphism, taking \(D_\lambda\) to
\(\Psi(\{(\lambda,\lambda)\})\) and \(D^c_\lambda\) to \(\Gg_\Lambda^{(0)}\cap\Psi(\{(\lambda,\lambda)\})^c\).
\end{rmk}

\begin{prop}
\label{prop:Psi is a hm} The map $\Psi$ defined by \eqref{eqn:Psi} is an injective
$*$-homomorphism of $S_{\Lambda}$ into $\Gg_{\Lambda}^{op}$.
\end{prop}
\begin{proof}
First we show that for $F,G\in S_{\Lambda}$, $\Psi(F)\Psi(G)=\Psi(FG)$.
Composition in $\Gg_{\Lambda}$ gives $\Psi(F)\Psi(G)\subseteq\Psi(FG)$.
For the other containment, let $[FG,x]\in\Psi(FG)$. Then there exists
$(\tau,\omega)\in FG$ such that $x(0,d(\omega))=\omega$. From the
definition of multiplication in $S_{\Lambda}$, there exist $(\lambda,\mu)\in F$,
$(\xi,\eta)\in G$ and $(\alpha,\beta)\in\Lmin(\mu,\xi)$ such that
$(\lambda\alpha,\eta\beta)=(\tau,\omega)$. But then $x(0,d(\eta\beta))=\eta\beta$,
so $x\in D_{G}$ and
\[
\theta_{G}(x)=\xi\bp(x;{d(\eta)},{d(x)})=\xi\beta\bp(x;{d(\eta\beta)},{d(x)})=\mu\alpha\bp(x;{d(\eta\beta)},{d(x)})\in
D_{F}.
\]
 Thus, $[FG,x]=[F,\theta_{G}(x)]\cdot[G,x]$, which gives $\Psi(FG)\subseteq\Psi(F)\Psi(G)$.

For each $F\in S_{\Lambda}$, we have $\Psi(F)\in\Gg_{\Lambda}^{op}$
since $\theta_{F}:D_{F}\rightarrow D_{F^{*}}$ is injective. Therefore
$\Psi$ is a homomorphism between inverse semigroups $S_{\Lambda}$
and $\Gg_{\Lambda}^{op}$, and so preserves involution by \cite[Proposition~2.1.1(iv)]{P1}.

Finally we show that $\Psi$ is injective. Fix $F,G\in S_{\Lambda}$
with $F\neq G$. Assume without loss of generality that there exists
$(\lambda,\mu)\in F$ such that $(\lambda,\mu)\not\in G$. Suppose
for contradiction that $\Psi(F)=\Psi(G)$. Then in particular, regarding
the path $\mu$ as an element of $X_{\Lambda}$, we have $\mu\in D_{F}=D_{G}$,
so there exists $(\xi,\eta)\in G$ such that $\mu(0,d(\eta))=\eta$
and \[
[\{(\lambda,\mu)\},\mu]=[F,\mu]=[G,\mu]=[\{(\xi,\eta)\},\mu].\]
 By \eqref{eqn:equivalence rel}, we then have \begin{equation}
\xi\mu(d(\eta),d(\mu))=\lambda.\label{eqn:xi=lambda}\end{equation}
 We claim that $(\lambda,\mu)=(\xi,\eta)$, which contradicts the
assumption that $(\lambda,\mu)\not\in G$.

Regarding $\eta$ as an element of $X_{\Lambda}$, we have $\eta\in D_{G}=D_{F}$,
so there exists $(\tau,\omega)\in F$ such that $\eta(0,d(\omega))=\omega$.
But then $\mu$ and $\omega$ have a common extension, namely \begin{equation}
\mu=\eta\mu(d(\eta),d(\mu))=\omega\eta(d(\omega),d(\eta))\mu(d(\eta),d(\mu)).\label{eqn:mu=eta}\end{equation}
 Hence by the definition of $S_{\Lambda}$ we must have $(\lambda,\mu)=(\tau,\omega)$,
from which \eqref{eqn:mu=eta} gives $\mu(d(\eta),d(\mu))=s(\omega)$.
Thus \eqref{eqn:xi=lambda} and \eqref{eqn:mu=eta} give $\lambda=\xi$
and $\mu=\eta$, contradicting our assumption that $(\lambda,\mu)\not\in G$.
\end{proof}
\begin{prop}
\label{prop:ample groupoid} $\Gg_{\Lambda}$ is an ample groupoid.
\end{prop}
\begin{proof}
For $G\in S_{\Lambda}$, $\Psi(G)^{c}=\bigcup_{F\in S_{\Lambda}}\Psi(F)\cap\Psi(G)^{c}$.
Hence $\{\Psi(F)\cap\Psi(G)^{c}:F,G\in S_{\Lambda}\}$ is a subbasis
for the topology of $\Gg_{\Lambda}$. Each $\Psi(F)$ is compact by
Proposition~\ref{prop:G is an S-groupoid}, and is an element of
$\Gg_{\Lambda}^{op}$ by Proposition~\ref{prop:Psi is a hm}. Hence
each $\Psi(F)\cap\Psi(G)^{c}$ is a compact element of $\Gg_{\Lambda}^{op}$.
The result follows.
\end{proof}
\begin{theorem}\label{thm:toeplitz isomorphism}
Let \((\Lambda,d)\) be a finitely aligned \(k\)-graph. Then the set of characteristic functions
\(\{1_{\Psi(\{(\lambda,s(\lambda))\})}  :  \lambda\in\Lambda\}\) is a Toeplitz-Cuntz-Krieger \(\Lambda\)-family in
\(C^*(\Gg_\Lambda)\) which gives a canonical isomorphism
\[
\Tt C^*(\Lambda)\cong C^*(\Gg_\Lambda).
\]
\end{theorem}

\begin{proof}
Proposition~\ref{prop:G is an S-groupoid} and \cite[Proposition~2.2.6]{P1} imply
that the map $\{(\lambda,\mu)\}\mapsto1_{\Psi(\{(\lambda,\mu)\})}$ defines a
$*$-homomorphism of $S_{\Lambda}$ into $C^{*}(\Gg_{\Lambda})$, and
Definition~\ref{C-K relations} $(1)-(3)$ follow from this.

By the universal property of $\Tt C^{*}(\Lambda)$ (Remark~\ref{rmk:universal properties}),
there is a homomorphism $\pi_{\Tt}:\Tt C^{*}(\Lambda)\rightarrow C^{*}(\Gg_{\Lambda})$
such that $\pi_{\Tt}(s_{\lambda}^{\Tt})=1_{\Psi(\{(\lambda,s(\lambda))\})}$.

Next we show that $\pi_\Tt$ is surjective; we do this in two steps, first showing that $C_c(\Gg_\Lambda^{(0)})$ is in
the image of $\pi_\Tt$.

Let $W = \lsp\{ 1_{\Psi(\{(\lambda,\lambda)\})} : \lambda\in\Lambda\} = \lsp\{\pi_\Tt(s_\lambda s^*_\lambda) :
\lambda\in \Lambda\}$. It is easy to see that $W$ is a $*$-subalgebra of $C_c(\Gg_\Lambda^{(0)})$. Furthermore, $W$
separates points in $\Gg_\Lambda^{(0)}$ and does not vanish identically at any point of $\Gg_\Lambda^{(0)}$. Therefore
by the Stone-Weierstrass Theorem, $W$ is uniformly dense in $C_c(\Gg_\Lambda^{(0)})$, and it follows that
$C_c(\Gg_\Lambda^{(0)})$ is in the image of $\pi_\Tt$.

Now fix $f\in C_c(\Gg_\Lambda)$. Let $\{\Psi(\{(\lambda_i,\mu_i)\})\}_{i=1}^n$ be a covering of $\supp f$ and let
$\{\phi_i\}_{i=1}^n$ be a partition of unity subordinate to $\{\Psi(\{(\lambda_i,\mu_i)\})\}_{i=1}^n$; that is,
$\phi_i:\Gg_\Lambda\to [0,1]$, $\supp\phi_i\subseteq \Psi(\{(\lambda_i,\mu_i)\})$ and $\sum_{i=1}^n\phi_i([F,x]) = 1$
for all $[F,x]\in \bigcup_{i=1}^n\Psi(\{(\lambda_i,\mu_i)\})$. We then have $f = \sum_{i=1}^n \phi_if$ (where the
product is pointwise) and $\supp (\phi_if) \subseteq \Psi(\{(\lambda_i,\mu_i)\})$; fixing $i\in\{1,\dots, n\}$, we will
show that $\phi_if$ is in the image of $\pi_\Tt$. Define $g\in C_c(\Gg_\Lambda^{(0)})$ by
\[
g(x) :=
\begin{cases}
(\phi_if)([\{(\lambda_i,\mu_i)\},\mu_ix]) &\text{if $r(x) = s(\mu_i)$} \\
0 &\text{otherwise.}
\end{cases}
\]
Then there exists $a\in \Tt C^*(\Lambda)$ such that $\pi_\Tt(a) = g$. Furthermore,
\[
\phi_i f = 1_{\Psi(\{(\lambda_i, s(\lambda_i))\})} \pi_\Tt(a) 1_{\Psi(\{(s(\mu_i),\mu_i)\})} =
\pi_\Tt(s_{\lambda_i})\pi_\Tt(a)\pi_\Tt(s_{\mu_i}^*),
\]
as required. Thus $C_c(\Gg_\Lambda)$ is in the image of $\pi_\Tt$. Since $\pi_\Tt$ is a $C^*$-homomorphism, and hence
closed, it follows that $\pi_\Tt$ is surjective.

To see that $\pi_{\Tt}$ is injective,
it suffices, by \cite[Theorem~8.1]{RS}, to show that for all $v\in\Lambda^{0}$
and all finite $E\subseteq(v\Lambda)\setminus\{ v\}$,
\begin{equation}
\prod_{\lambda\in E}\big(1_{\Psi(\{(v,v)\})}-1_{\Psi(\{(\lambda,s(\lambda))\})}
1_{\Psi(\{(\lambda,s(\lambda))\})}^{*}\big)>0,\label{eqn:TCK uniqueness rel}
\end{equation}
 where the product is convolution. For $v\in\Lambda^{0}$ and finite
$E\subseteq(v\Lambda)\setminus\{ v\}$, regarding $v$ as an element
of $X_{\Lambda}$, we have $[\{(v,v)\},v]\in\Psi(\{(v,v)\})$ but
$[\{(v,v)\},v]\not\in\Psi(\{(\lambda,\lambda)\})$ for all $\lambda\in E$.
Thus the product in \eqref{eqn:TCK uniqueness rel} is bounded below
by $1_{\{[\{(v,v)\},v]\}}$ which is nonzero in $C^{*}(\Gg_{\Lambda})$
since the left Haar system on $\Gg_{\Lambda}$ is given by the counting
measures.
\end{proof}
\begin{prop}
\label{lem:pL closed invariant} $\partial\Lambda$ is a nonempty closed invariant
subset of $\Gg_{\Lambda}^{(0)}$.
\end{prop}
\begin{proof}
Since $\Gg_{\Lambda}^{(0)}$ is homeomorphic to $X_{\Lambda}$ (see
Remark~\ref{rmk:unit space topology}) and since $r([F,x])=\theta_{F}(x)$,
Remark~\ref{rmk:partial Lambda closed invariant} implies that $\partial\Lambda$
is a closed invariant subset of $\Gg_{\Lambda}^{(0)}$. It is nonempty by Lemma \ref{lem:v partialLambda nonempty}.
\end{proof}
Proposition~\ref{lem:pL closed invariant} implies that $\Gg_{\Lambda}|_{\partial\Lambda}$
is a locally compact, Hausdorff and ample groupoid.

For $F\in S_{\Lambda}$, define $\Psi_{*}:S_{\Lambda}\rightarrow(\Gg_{\Lambda}|_{\partial\Lambda})^{a}$
by \begin{equation}
\Psi_{*}(F):=\Psi(F)\cap\Gg_{\Lambda}|_{\partial\Lambda}=\{[F,x]:x\in
D_{F}\cap\partial\Lambda\}\label{eqn:Psi_*}\end{equation}

\begin{rmks}\label{rmk:restricted groupoid notes}
\begin{itemize}
\item[(1)]
Unlike $\Psi$, \(\Psi_*\) is not injective in general. For example,
for the 1-graph given by
\[
\begin{picture}(60,10)
\put(0,0){\(\bullet\)}
\put(24,2.5){\(\vector(-1,0){20}\)}
\put(23,0){\(\bullet\)}
\put(47,2.5){\(\vector(-1,0){20}\)}
\put(46,0){\(\bullet\)}
\put(53,0){\(\cdots\)}
\put(0,5){\(v\)}
\put(12,5){\(e_1\)}
\put(34,5){\(e_2\)}
\end{picture}
\]
we have \(\Psi_*(\{(v,v)\})=\Psi_*(\{(e_1,e_1)\})=\Psi_*(\{(e_1\cdots e_n,e_1\cdots e_n)\})\).
In general, given a finitely aligned \(k\)-graph \(\Lambda\), for all \(v\in\Lambda^0\) and \(E\in v\FE(\Lambda)\),
\[
\Psi_*(\{(v,v)\})=\bigcup_{\lambda\in E}\Psi_*(\{(\lambda,\lambda)\}).
\]
\item[(2)]
The family of sets of the form
\[
\Psi_*(\{(\lambda,\mu)\})\cap \Psi_*(\{(\lambda\nu_1,\mu\nu_1)\})^c\cap\cdots\cap\Psi_*(\{(\lambda\nu_l,\mu\nu_l)\})^c
\]
where \((\lambda,\mu)\in\Lambda *_s\Lambda\) and \(\nu_1,\dots,\nu_l\in s(\lambda)\Lambda\), is a basis for the
topology on \(\Gg_\Lambda|_{\partial\Lambda}\).
\item[(3)]
In \cite{KP}, the authors consider the \(C^*\)-algebras arising from \(k\)-graphs which are row-finite and have no
sources (Recall that this means \(0<|v\Lambda^n|<\infty\) for all \(v\in\Lambda^0\) and $n\in\NN^k$.) by constructing
and analyzing a groupoid we call \({\mathcal G}^{KP}_\Lambda\) (see \cite[Definition~2.7]{KP}). In this setting, the
boundary paths \(\partial\Lambda\) are precisely the paths \(x\in X_\Lambda\) such that \(d(x)=(\infty,\dots,\infty)\)
(that is, \(x\) is infinite in all \(k\) directions), and using the description of \(\Gg_\Lambda\) given in
Remark~\ref{rmk:composable pairs} and Remark~\ref{rmk:restricted groupoid notes}~(2), one can see that the reduction
\(\Gg_\Lambda|_{\partial\Lambda}\) is isomorphic as a topological groupoid to \({\mathcal G}^{KP}_\Lambda\).
\end{itemize}
\end{rmks}

For $z\in\TT^{k}$ and $m\in\NN^{k}$, define $z^{m}:=z_{1}^{m_{1}}z_{2}^{m_{2}}\cdots z_{k}^{m_{k}}\in\TT$.
There is a strongly continuous \emph{gauge action} $\gamma$ of $\TT^{k}$
on $C^{*}(\Lambda)$ defined by $\gamma_{z}(s_{\lambda}s_{\mu}^{*})=z^{d(\lambda)-d(\mu)}s_{\lambda}s_{\mu}^{*}$
for $z\in\TT^{k}$. The gauge-invariant uniqueness theorem for finitely
aligned $k$-graphs \cite[Corollary~4.3]{RSY2} says that if $A$
is a $C^{*}$-algebra generated by a Cuntz-Krieger $\Lambda$-family,
and if $A$ carries a strongly continuous action of $\TT^{k}$ which
is equivariant to the gauge action on $C^{*}(\Lambda)$, then $A$
is isomorphic to $C^{*}(\Lambda)$. The next lemma shows that $C^{*}(\Gg_{\Lambda}|_{\partial\Lambda})$
admits such an action.

\begin{lemma}\label{lem:action on Gg_Lambda}
There is a strongly continuous action \(\beta:\TT^k\to \Aut(C^*(\Gg_\Lambda|_{\partial\Lambda}))\) such that
\(\beta_z(1_{\Psi_*(\{(\lambda,s(\lambda))\})})=z^{d(\lambda)}1_{\Psi_*(\{(\lambda,s(\lambda))\})}\) for all
\(\lambda\in\Lambda\).
\end{lemma}

\begin{proof}
Define a map $c:\Gg_{\Lambda}|_{\partial\Lambda}\rightarrow\ZZ^{k}$
by $[\{(\lambda,\mu)\},x]\mapsto d(\lambda)-d(\mu)$. To see that
$c$ is well-defined, suppose $[\{(\lambda,\mu)\},x]=[\{(\xi,\eta)\},x]$.
Then by \eqref{eqn:equivalence rel} we have \[
\lambda x(d(\mu),d(\mu)\vee d(\eta))=\xi x(d(\eta),d(\mu)\vee d(\eta)),\]
 so \begin{align*}
d(\lambda)-d(\mu) & =d(\lambda)+(d(\mu)\vee d(\eta))-d(\mu)-(d(\mu)\vee d(\eta))\\
 & =d(\xi)+(d(\mu)\vee d(\eta))-d(\eta)-(d(\mu)\vee d(\eta))\\
 & =d(\xi)-d(\eta).\end{align*}
 To see that $c$ is a 1-cocycle, take $\left(\left[\{(\lambda,\mu)\},\xi\bp(x;{d(\eta)},{d(x)})\right],
\left[\{(\xi,\eta)\},x\right]\right)\in\Gg_{\Lambda}|_{\partial\Lambda}^{(2)}$.
Let
\[
\alpha:=x(d(\eta),d(\eta)+(d(\mu)\vee d(\xi))-d(\mu))=\big(\bp(x;{d(\eta)},{d(x)})\big)(0,(d(\mu)\vee d(\xi))-d(\mu))
\]
and
\[
\beta:=x(d(\eta),d(\eta)+(d(\mu)\vee d(\xi))-d(\xi))=\big(\bp(x;{d(\eta)},{d(x)})\big)(0,(d(\mu)\vee d(\xi))-d(\xi)).
\]
Then $(\alpha,\beta)\in\Lmin(\mu,\xi)$, and we have
\begin{align*}
c & \big(\big[\{(\lambda,\mu)\},\xi\bp(x;{d(\eta)},{d(x)})\big]\cdot[\{(\xi,\eta)\},x]\big)\\
 & =c([\{(\lambda\alpha,\eta\beta)\},x])\\
 & =d(\lambda\alpha)-d(\eta\beta)\\
 & =\big(d(\lambda)+(d(\mu)\vee d(\xi))-d(\mu)\big)-\big(d(\eta)+(d(\mu)\vee d(\xi))-d(\xi)\big)\\
 & =d(\lambda)-d(\mu)+d(\xi)-d(\eta)\\
 & =c\big(\big[\{(\lambda,\mu)\},\xi\bp(x;{d(\eta)},{d(x)})\big]\big)+c([\{(\xi,\eta)\},x]).\end{align*}
 Hence $c$ is a $\ZZ^{k}$-valued 1-cocycle. To see that $c$ is
continuous, simply observe that for $n\in\ZZ^{k}$ \[
c^{-1}(\{ n\})=\bigcup_{\substack{\{(\lambda,\mu)\}\in S_{\Lambda}\\
d(\lambda)-d(\mu)=n}
}\Psi_{*}(\{(\lambda,\mu)\})\]
 is open in $\Gg_{\Lambda}|_{\partial\Lambda}$. By \cite[II.5.1]{R},
there is a strongly continuous action $\beta:\TT^{k}\rightarrow\Aut(C^{*}(\Gg_{\Lambda}|_{\partial\Lambda}))$
such that \[
\beta_{z}(1_{[\{(\lambda,\mu)\},x]})=z^{d(\lambda)-d(\mu)}1_{[\{(\lambda,\mu)\},x]}\]
 for all $[\{(\lambda,\mu)\},x]\in\Gg_{\Lambda}|_{\partial\Lambda}$.
The result follows.
\end{proof}
\begin{theorem}\label{thm:CK isomorphism}
Let \((\Lambda,d)\) be a finitely aligned \(k\)-graph. Then the set of characteristic functions
\(\{1_{\Psi_*(\{(\lambda,s(\lambda))\})}  :  \lambda\in\Lambda\}\) is a Cuntz-Krieger \(\Lambda\)-family in
\(C^*(\Gg_\Lambda|_{\partial\Lambda})\) that gives a canonical isomorphism
\[
C^*(\Lambda)\cong C^*(\Gg_\Lambda|_{\partial\Lambda}).
\]
\end{theorem}

\begin{proof}
Proposition~\ref{prop:G is an S-groupoid} and \cite[Proposition~2.2.6]{P1} imply
that $\{(\lambda,\mu)\}\mapsto1_{\Psi_{*}(\{(\lambda,\mu)\})}$ defines a
$*$-homomorphism of $S_{\Lambda}$ into $C^{*}(\Gg_{\Lambda}|_{\partial\Lambda})$,
and (1)--(3) of Definition~\ref{C-K relations} follow from this. To show
Definition~\ref{C-K relations} (CK) holds, we fix $v\in\Lambda^{0}$, $E\in
v\FE(\Lambda)$ and $[F,x]\in\Gg_{\Lambda}|_{\partial\Lambda}$, and evaluate
\begin{equation} \left(\prod_{\lambda\in
E}\big(1_{\Psi_{*}(\{(v,v)\})}-1_{\Psi_{*}(\{(\lambda,\lambda)\})}\big)\right)[F,x]\label{eqn:product}\end{equation}
 where the product is convolution.

For $\lambda\in E$, $\Psi_{*}(\{(\lambda,\lambda)\})$ is a subset
of $\Psi_{*}(\{(v,v)\})$, so \eqref{eqn:product} can only be nonzero
if $[F,x]=[\{(v,v)\},x]$, in which case $x\in v(\partial\Lambda)$.
Since \eqref{eqn:product} is a product of characteristic functions
of subsets of $(\Gg_{\Lambda}|_{\partial\Lambda})^{(0)}$, a straightforward
calculation shows that \eqref{eqn:product} can be written as the
pointwise product \begin{equation}
\prod_{\lambda\in
E}\left(\big(1_{\Psi_{*}(\{(v,v)\})}-1_{\Psi_{*}(\{(\lambda,\lambda)\})}\big)[\{(v,v)\},x]\right).\label{eqn:pointwise
product}\end{equation}
 Since $E\in v\FE(\Lambda)$, there exists $\mu\in E$ such that $x(0,d(\mu))=\mu$.
We then have $[\{(v,v)\},x]=[\{(\mu,\mu)\},x]$, so the term
$\big(1_{\Psi_{*}(\{(v,v)\})}-1_{\Psi(\{(\mu,\mu)\})}\big)[\{(v,v)\},x]$
from \eqref{eqn:pointwise product} is zero. Therefore \eqref{eqn:product}
is zero, and $\{1_{\Psi_{*}(\{(\lambda,s(\lambda)\})}:\lambda\in\Lambda\}$
is a Cuntz-Krieger $\Lambda$-family. By the universal property of
$C^{*}(\Lambda)$ (Remark~\ref{rmk:universal properties}), there
is a unique homomorphism $\pi:C^{*}(\lambda)\rightarrow C^{*}(\Gg_{\Lambda}|_{\partial\Lambda})$
such that $\pi(s_{\lambda})=1_{\Psi_{*}(\{(\lambda,s(\lambda)\})}$
for all $\lambda\in\Lambda$. Furthermore, since each $v(\partial\Lambda)$
is nonempty by Proposition~\ref{lem:v partialLambda nonempty}, and
since the left Haar system on $\Gg_{\Lambda}|_{\partial\Lambda}$
is given by counting measures, it follows that $\pi(s_{v})=1_{\Psi_{*}(\{(v,v)\})}$
is nonzero for every $v\in\Lambda^{0}$.

Lemma~\ref{lem:action on Gg_Lambda} gives a strongly continuous
action $\beta:\TT^{k}\rightarrow\Aut(C^{*}(\Gg_{\Lambda}|_{\partial\Lambda}))$
such that $\beta_{z}\circ\pi=\pi\circ\gamma_{z}$ for all $z\in\TT^{k}$.
Therefore \cite[Theorem~4.2]{RSY2} implies that $\pi$ is injective. We can see that $\pi$ is surjective in the same way
that we saw that $\pi_\Tt$ was surjective in the proof of Theorem~\ref{thm:toeplitz isomorphism}. The result follows.
\end{proof}

\section{The Cuntz-Krieger Uniqueness Theorem}

\label{sec:CK uniqueness}

In this section we present a Cuntz-Krieger uniqueness theorem for
finitely aligned $k$-graphs. The theorem generalizes \cite[Theorem~4.6]{KP}
but differs from the Cuntz-Krieger uniqueness theorem \cite[Theorem~4.5]{RSY2}
obtained by direct methods. These differences will be discussed in
Remarks~\ref{rmk:CKUT differences}.

\begin{theorem}\label{thm:groupoid CK uniqueness}
Let \((\Lambda,d)\) be a finitely aligned \(k\)-graph, and suppose that
\begin{equation}\label{eqn:aperiodic}
\begin{split}
\text{for each } v\in\Lambda^0& \text{ there exists \(x\in v(\partial\Lambda)\) such that}\\
&\bp(x;m,{d(x)})=\bp(x;n,{d(x)}) \implies m=n \text{ for all \(m,n\in\NN^k\), \(m,n\le d(x)\).}
\end{split}\tag{A}
\end{equation}
Suppose that \(\pi:C^*(\Lambda)\to B\) is a homomorphism such that \(\pi(s_v)\) is nonzero for all \(v\in\Lambda^0\).
Then \(\pi\) is injective.
\end{theorem}

The proof of Theorem~\ref{thm:groupoid CK uniqueness} relies on
the following proposition.

\begin{prop}
\label{prop:ess free iff A} Let $(\Lambda,d)$ be a finitely aligned
$k$-graph. Then $\Lambda$ satisfies condition (A) if and only if
$\Gg_{\Lambda}|_{\partial\Lambda}$ is essentially free.
\end{prop}
\begin{proof}
First, observe that $x\in(\Gg_{\Lambda}|_{\partial\Lambda})^{(0)}$
has trivial isotropy if and only if \begin{equation}
\text{for all }m,n\in\NN^{k}\text{ with }m,n\leq d(x),\sigma^{m}x=\sigma^{n}x\text{ implies }m=n.\label{eqn:apperiodic
condition}\end{equation}

If $\Gg_{\Lambda}|_{\partial\Lambda}$ is essentially free, then since
each $D_{v}$ is nonempty by Lemma~\ref{lem:v partialLambda nonempty}
it follows that each $D_{v}=v(\partial\Lambda)$ contains a boundary
path $x$ satisfying \eqref{eqn:apperiodic condition}. Therefore
$\Lambda$ satisfies Condition~\eqref{eqn:aperiodic}.

Conversely, suppose $\Lambda$ satisfies Condition~\eqref{eqn:aperiodic},
and let $x\in D_{\lambda}\cap D_{\lambda\nu_{1}}^{c}\cap\cdots\cap
D_{\lambda\nu_{l}}^{c}\subseteq(\Gg_{\Lambda}|_{\partial\Lambda})^{(0)}$.

Since $x$ is a boundary path, for any finite exhaustive set $E\subseteq s(\lambda)\Lambda$
there exists $\xi\in E$ such that $x(d(\lambda),d(\lambda)+d(\xi))=\xi$.
Hence $\{\nu_{1},\dots,\nu_{l}\}$ cannot be exhaustive, so there
exists $\eta\in s(\lambda)\Lambda$ such that $\Lmin(\eta,\nu_{j})=\emptyset$
for all $j\in\{1,\dots,l\}$. Furthermore, since $\Lambda$ satisfies
Condition~\eqref{eqn:aperiodic}, there exists $y\in s(\eta)(\partial\Lambda)$
satisfying \[
\sigma^{m}y=\sigma^{n}y\text{ implies }m=n\text{ for all }m,n\leq d(y).\]
 Therefore $\lambda\eta y$ is an element of $D_{\lambda}\cap D_{\lambda\nu_{1}}^{c}\cap\cdots\cap
D_{\lambda\nu_{l}}^{c}$
satisfying \[
\sigma^{m}(\lambda\eta y)=\sigma^{n}(\lambda\eta y)\text{ implies }m=n\text{ for all }m,n\leq d(\lambda\eta y).\]
 Since $D_{\lambda}\cap D_{\lambda\nu_{1}}^{c}\cap\cdots\cap D_{\lambda\nu_{l}}^{c}$
was an arbitrary basis set containing $x$, it follows that $\Gg_{\Lambda}|_{\partial\Lambda}$
is essentially free.
\end{proof}

\begin{proof}[{Proof of Theorem~\ref{thm:groupoid CK uniqueness}}]
Let $\pi_*: C^*(\Lambda)\to C^*(\Gg_\Lambda|_{\partial\Lambda})$ be the
canonical isomorphism of Theorem~\ref{thm:CK isomorphism}.
Then $\pi$ is injective if and only if $\pi\circ \pi_*^{-1}:
C^*(\Gg_\Lambda|_{\partial\Lambda})\to B$ is injective. By
Proposition~\ref{prop:ess free iff A},
$\Gg_\Lambda|_{\partial\Lambda}$ is essentially free. Moreover,
$\Gg_\Lambda|_{\partial\Lambda}$ is $r$-discrete and
$(\Gg_\Lambda|_{\partial\Lambda})^{(0)}$ has a basis consisting of
compact open sets, so by \cite[Corollary~3.6]{KPR} it suffices to show
that $\pi\circ \pi_*^{-1}$ is injective on
$C_0((\Gg_\Lambda|_{\partial\Lambda})^{(0)})$. If the kernel of the
restriction of $\pi\circ\pi_*^{-1}$ to
$C_0((\Gg_\Lambda|_{\partial\Lambda})^{(0)})$ is nonzero, it must
contain a characteristic function
$1_{D_\lambda}=1_{\Psi_*(\{(\lambda,\lambda)\})}$ for some
$\lambda\in\Lambda$. It follows that $\pi(s_\lambda s^*_\lambda)=0$,
which implies $\pi(s^*_\lambda s_\lambda) = \pi(s_{s(\lambda)}) = 0$, a
contradiction.

\end{proof}

\begin{rmks}\label{rmk:CKUT differences}
There are a number of issues to point out here. Firstly, when the \(k\)-graph is row-finite and has no sources, our
condition \eqref{eqn:aperiodic} is equivalent to the aperiodicity condition of \cite[Definition~4.3]{KP}, and
Theorem~\ref{thm:groupoid CK uniqueness} gives \cite[Theorem~4.6]{KP}.

In \cite{RSY2}, a Cuntz-Krieger uniqueness theorem is given for finitely aligned
\(k\)-graphs \cite[Theorem~4.5]{RSY2}. To compare \cite[Theorem~4.5]{RSY2} with
Theorem~\ref{thm:groupoid CK uniqueness}, define \(\Lambda^{\le\infty}\) to be the
subset of \(X_\Lambda\) consisting of all \(x\in X_\Lambda\) for which there exists
\(n_x\in\NN^k\), \(n_x\le d(x)\), satisfying
\[
\begin{split}
n\in\NN^k, n_x\le n\le d(x) \text{ and } n_i=d(x)_i \text{ imply that } x(n)\Lambda^{e_i}=\emptyset,
\end{split}
\]
and for \(v\in\Lambda^0\), define \(v\Lambda^{\le\infty}:=\{x\in\Lambda^{\le\infty}
: r(x)=v\}\). Then \cite[Theorem~4.5]{RSY2} differs from Theorem~\ref{thm:groupoid
CK uniqueness} in that our condition~\eqref{eqn:aperiodic} is replaced with the
condition
\begin{equation}\label{eqn:RSY apperiodicity}
\begin{split}
\text{for each }& v\in\Lambda^0  \text{ there exists } x\in v\Lambda^{\le\infty} \text{ such that} \\
\lambda,\mu&\in\Lambda v \text{ and } \lambda\neq\mu \text{ implies } \lambda x\neq \mu x.
\end{split}\tag{B}
\end{equation}
The two conditions \eqref{eqn:aperiodic} and \eqref{eqn:RSY apperiodicity} do not
seem equivalent \emph{a priori}; when the \(k\)-graph \(\Lambda\) is row-finite and
has no sources, the set \(\Lambda^{\le\infty}\) is precisely \(\partial\Lambda\),
and condition~\eqref{eqn:aperiodic} implies condition~\eqref{eqn:RSY apperiodicity} (see \cite[Remark~4.4]{RSY1}),
making \eqref{eqn:RSY apperiodicity} seem the weaker
condition. However, even in the row-finite and no sources setting, it remains
unclear whether \eqref{eqn:RSY apperiodicity} is strictly weaker than
\eqref{eqn:aperiodic}.

When the \(k\)-graph \(\Lambda\) is finitely aligned, it is not clear that either
condition implies the other. The set \(\Lambda^{\le\infty}\) used in \cite{RSY2} is
in general a proper subset of \(\partial\Lambda\). If \(\Lambda^{\le\infty}\) is
replaced with the set \(\partial\Lambda\) in \eqref{eqn:RSY apperiodicity} to give
\begin{equation}\label{eqn:N}
\begin{split}
\text{for each }& v\in\Lambda^0  \text{ there exists } x\in v(\partial\Lambda) \text{ such that} \\
\lambda,\mu&\in\Lambda v \text{ and } \lambda\neq\mu \text{ implies } \lambda x\neq \mu x,
\end{split}\tag{B\protect\('\protect\)}
\end{equation}
then the resulting property of the groupoid \(\Gg_\Lambda|_{\partial\Lambda}\) is
not the essential freeness implied by \eqref{eqn:aperiodic}, rather the curious
property:
\begin{gather*}
\text{for all } v\in\Lambda^0, \text{ there exists } x\in v(\partial\Lambda) \text{ such that } \\
r([\{(\lambda,v)\},x]) = r([\{(\mu,v)\},x]) \text{ implies } [\{(\lambda,v)\},x]=[\{(\mu,v)\},x].
\end{gather*}
Although this property of \(\Gg_\Lambda|_{\partial\Lambda}\) is quite different to
essential freeness, it may still yield similar consequences to essential freeness.
In particular, the conclusion \cite[Lemma~3.5]{KPR} may hold, the key result in the
proof of \cite[Theorem~4.6]{KP} and Theorem~\ref{thm:groupoid CK uniqueness}. If
this were the case, and if condition \eqref{eqn:RSY apperiodicity} was not
equivalent to \eqref{eqn:N}, then we would obtain a generalization
of \cite[Theorem~4.5]{RSY2} using groupoid methods (since \eqref{eqn:RSY apperiodicity} asks for an element of
\(v(\partial\Lambda)\)
for each \(v\in\Lambda^0\), whereas \eqref{eqn:N} asks for an element of the smaller
set \(v\Lambda^{\le\infty}\)).
\end{rmks}

\end{document}